\documentclass[twocolumn, showpacs, preprintnumbers, amsmath, amssymb]{revtex4}
\usepackage{pifont}
\usepackage{color}
\usepackage{graphicx}% Include figure files
\usepackage{dcolumn}% Align table columns on decimal point
\usepackage{bm}% bold math
\usepackage{appendix}
\begin{document}

\preprint{APS/123-QED}

\title{Effects of the tempered aging and its Fokker-Planck equation}

\author{Weihua Deng}
\author{Wanli Wang}
\author{Xinchun Tian}
\author{Yujiang Wu}
 %\altaffiliation[Also at ]{School of Mathematics and Statistics,  Lanzhou University. }%Lines break automatically or can be forced with \\
%\email{Second. Author@institution. edu}

% \author{Eli Barkai$^2$}%

\affiliation{%
School of Mathematics and Statistics, Gansu Key Laboratory of Applied Mathematics and Complex
Systems, Lanzhou University,
Lanzhou 730000,  P.R. China
%\\
%$^2$Department of Physics, Advanced Materials and Nanotechnology Institute, Bar Ilan University,  Ramat-Gan 52900,
%Israel
 }%

% \homepage{http://www. Second. institution. edu/~Charlie. Author}
%\affiliation{

%\date{\today}% It is always \today,  today,

\begin{abstract}
In the renewal processes, if the waiting time probability density function is a tempered power-law distribution, then the process displays a transition dynamics; and the transition time depends on the parameter $\lambda$ of the exponential cutoff. In this paper, we discuss the aging effects of the renewal process with the tempered power-law waiting time distribution. By using the aging renewal theory, the $p$-th moment of the number of renewal events $n_a(t_a, t)$ in the interval $(t_a, t_a+t)$ is obtained for both the weakly and strongly aged systems; and the corresponding surviving probabilities are also investigated. We then further analyze the tempered aging continuous time random walk and its Einstein relation, and the mean square displacement is attained. Moreover, the tempered aging diffusion equation is derived.

%In this paper we use the exponentially tempered waiting time possibility density function (PDF) to describe the aging effects.  we use the aging renewal theory to obtain the pth moment of the renewal process $n_a(t_a, t)$ in the interval $(t_a, t_a+t)$ and analyze two different systems; namely; the slightly  and the strongly aged systems. We then study  tempered aging continuous random walks~(ACTRW) and get the mean square displacement  for the two limit situations.  Moreover we discuss the aging diffusion equation with tempered; which is a fractional diffusion equation and can describe more generally tempered ACTRW process.
\end{abstract}

% PACS,  the Physics and Astronomy
\maketitle

\section{Introduction}\label{sec1}
In 1975, Scher and Motroll \cite{Scher:1} used the continuous time random walk (CTRW) to study non-Gaussian anomalous diffusion. Nowadays, the CTRW model becomes popular in describing anomalous diffusion and a lot of chemical,  physical, and biological processes \cite{Naftaly:1,Barkai:2,Barkai:6}, such as, the transport of electric charge in a complex system, diffusion in a low dimensional chaotic system, and the anomalous diffusion when cooling the  mental solid and the twinkling of single quantum dot. In 1996, Monthusyx and Bouchaud introduce a CTRW framework for describing the aging phenomena in glasses \cite{Monthusyx:1}.  This generalized CTRW is called aging continuous time random walk (ACTRW) in \cite{Barkai:1}. The complex dynamical systems displaying aging behaviour are quite extensive, including the fluorescence of single nanocrystals \cite{Brokmann:1}, aging effect in a single-particle trajectory averages \cite{Schulz:2}.

%Aging is a complex process. All kinds of organisms grow with time, finally lose their functions, and even die. Some of the organisms display anomalous diffusion \cite{Cherstvya:2,Klafter:2}.

%Aging is ubiquitous in biology, psychology, medical science, and so on. In our daily life, it can be easily noted a large number of aging phenomena, e.g., the aging of people, bulb, and equipment, etc. In fact,

% Aging is useful and effective to research the complex system, such as, spinglass, colloidal, the twinkling of single quantum dot. Today, a list of complex dynamics systems displaying aging behaviour is quite extensive: the fluorescence of single nanocrystals\cite{Brokmann:1}, aging effect in a single-particle trajectory averages\cite{Schulz:2}.

Research on statistics, based on  power-law distributions with a heavy tail, yields
many of significant results. Often the power-law distribution doesn't extend indefinitely, due to the finite life span of particles, the boundedness of physical space. For this reason, in 1994, Mantegna and Stanley  omit the large steps to study the truncated L\'{e}vy flights \cite{Mantegna:1,Negrete:1}. While the tempered power-law distribution \cite{Rosinski:1} uses a different approach, exponentially tempering the probability of large jumps.  Exponential tempering offers technical advantages since the tempered process is still an infinitely divisible L\'{e}vy process which makes it convenient to identify the governing equation and compute the transition densities at any scale \cite{Baeumera:1}. By tempering, the distribution changes from heavy tail to semi-heavy, and the existence of conventional moments is ensured, which is useful in some practical applications. Recently tempered power-law distributions \cite{Meerschaert:1} have been observed for many geophysical processes at various scales \cite{Baeumera:1,Negrete:1,Sokolov:1,Meerschaert:2,Allegrini:1}, including interplanetary solar-wind velocity and magnetic field fluctuations measured in the alluvial aquifers \cite{Bruno:1}.

%offers technical advantages, since the tempered process remains an infinitely divisible L¨¦vy process whose governing equation can be identified, and whose transition densities can be computed at any scale; those transition densities solve a tempered fractional diffusion equation, quite similar to the fractional diffusion equation [15]

%Previous research on statistical  which is based on the power-law distributed with a heavy tail waiting time comes out many of  significative results. while often the power law distribution doesn't extend indefinitely , due to the standard central limit doesn't hold and the infinite variance. for this reason, In 1994, R. N. Mantegna and H. E.  Stanley simply omit the large steps to study the truncated L$\acute{e}$vy flights \cite{Mantegna:1}, while the tempered power distribution use a different approach, the large steps are cut down by the  exponential distribution. The distribution tail change from heavy tail to semi-heavy, and the existence of conventional moments is ensured, which is useful in practical application. Recently tempered power-law distribution have been observed for many geophysical processes at various scales\cite{Baeumera:1,Negrete:1,Sokolov:1,Meerschaert:2,Allegrini:1}, including interplanetary solar-wind velocity and magneticfield fluctuations measured in the alluvial aquifers. We choose the PDF Eq.(\ref{ageq11}), which can capture the natural cutoff of  times, especially, which serves as a smoother alternative without a sharp cutoff.

In this paper, we discuss the aging effects of the renewal processes with exponentially tempered power-law waiting time  probability density function (PDF)
\begin{equation}\label{ageq11}
\varphi(t)= L_\alpha(t) e^{\lambda^\alpha-\lambda t} \sim \frac{1}{-\Gamma(-\alpha)}t^{-(1+\alpha)}e^{-\lambda t},
\end{equation}
 where $0<\alpha<1$, $L_\alpha(t)$ is the one side L\'{e}vy distribution \cite{Feller:1,levy:1}, and $\lambda>0$ is generally a small parameter.  The semi-heavy tails and scale-free waiting time properties of $\varphi(t)$ play a particularly prominent role in  diffusion phenomena.

\begin{figure}[htb]
  \centering
  \includegraphics[width=9cm, height=6cm]{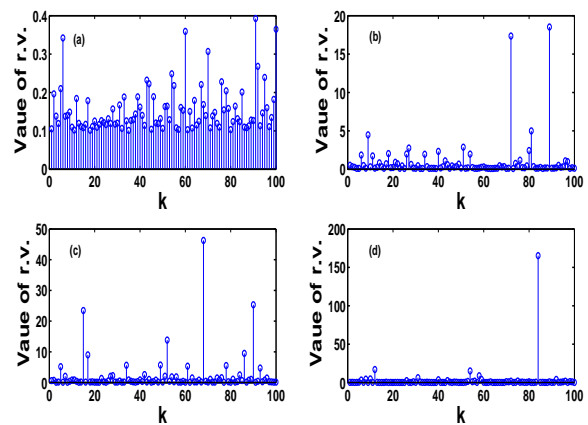}\\
  \caption{Random variables (r.v.) generated by Eq. (\ref{agfig11}) with $\alpha=0.6$. And the parameter $\lambda$ is chosen, respectively, as (a) $\lambda=10$, (b) $\lambda=10^{-1}$, (c) $\lambda=10^{-3}$, (d) $\lambda=10^{-5}$ and $k=1$, $2$, $\cdots$, 100.
}\label{agfig11}
\end{figure}

From Eq. (\ref{ageq11}), it can be noted that if $0.1\ll t\ll 1/\lambda$, $\varphi(t)\sim t^{-(1+\alpha)}$, while if $t\gg 1/\lambda$, $\varphi(t)\sim \exp(-\lambda t)$. For the random variables generated by Eq. (\ref{ageq11}), Fig. \ref{agfig11} shows that the maximum and range of fluctuations vary dramatically with the change of $\lambda$. The introduced tempering forces the renewal process to converge from non-Gaussian to Gaussian. But the convergence is very slow, requiring a long time to find the trend. So, with the time passed by, both the non-Gaussian and Gaussian processes can be described.

  %From Eq.(\ref{ageq11}), we have two limits form, when $t\ll 1/\lambda$ we have $\varphi(t)\sim t^{-(1+\alpha)}$, with time pass by, when $t\gg 1/\lambda$, $exp(-\lambda t)$ paly a important role. and the introduce of a tempered force the truncated diffusion change from anomalous diffusion to normal diffusion, thus get rid of the divergent variance. In addition, the convergence is very slow, requiring a long to time to find the trend. so we can describe  both the non-Gaussian and Gaussian process with time passed by.

  %Numerical simulation of $t_{i}$ defined by Eq. (\ref{agfig11}), the index of $\lambda$ is or  given by (a) $\lambda=10$, (b) $\lambda=1$, (c) $\lambda=0. 01$, (d) $\lambda=0. 0001$.  For large $\lambda$ ($\lambda>0. 1$) the size of step length is smaller than $0. 35$, the length of each steps makes litter difference. while for small $\lambda=10^{-4}$, the large step length is about $10$ and the length of steps vary dramatically.  Different $\lambda$ has different character.

This paper is organized as follows. In Sec. \ref{sec2}, we use the aging renewal theory to obtain the $p$-th moment of the number of renewals $n_a(t_a, t)$ within the time interval $(t_{a},t+t_{a})$ of the renewal process starting from time zero. The survival probability is, respectively, discussed in weakly and strongly aging system. We then turn to discuss the tempered ACTRW in Sec. \ref{sec3}, and the mean square displacement is obtained for the cases $t_a \ll t$ and $t \ll t_a$. The numerical simulations confirm the analytical expressions of the mean square displacement. And the propagator function is also numerically obtained. In Sec. \ref{sec4}, we discuss the Einstein relation of the tempered ACTRW. The diffusion equation for the tempered ACTRW is derived in Sec. \ref{sec5}, describing the time evolution of the PDF of the position. Finally, we conclude the paper with some remarks.

%This paper is organized as follows, In Sec. (\ref{sec2}), we use the aging renewal theory to obtain the $p$-th moment of the number of renewal $n_a(t_a, t)$ with in the time interval $(t_{a},t+t_{a})$ after the process starts at time zeros, and discuss the properties of $\langle n_a(t_a, t)\rangle$. We analyze the survival probability in weekly and strongly aged systems, respectively. We then turn to describe the tempered ACTRW in Sec.(\ref{sec3}). In this case of waiting time distribution, we obtain the mean square displacement in the two limits, and we can demonstrate effectiveness of the approach. In Sec.(\ref{sec4}) we discuss the propagator by the numerical inversion of Laplace transform. In Sec.(\ref{sec6}) we consider the Einstein relation for the tempered. The aging diffusion equation was studied in Sec.(\ref{sec5}), this is a fractional diffusion equation and can describe a large class of aging phenomenon. Finally, we briefly present the conclusion of the paper.

\section{Tempered aging renewal theory}\label{sec2}

First, we briefly outline the main ingredients in the CTRW and ACTRW. The standard CTRW assumes that the jumping transitions begin at time $t=0$, and  observation of dynamics starts at $t_{a}$. The ACTRW modifies the statistic of time interval for first jump, namely, the waiting time PDF to the first jump is $\omega(t_a,t)$.  It describes a CTRW process having the aging time interval $(0,t_a)$,  while $t_a$ corresponding to the initial observation time $t=0$. Aging means that the number of renewals in the  time interval $(t_{a}, t_{a}+t)$  depends on the aging time $t_a$, even when the former is long. Thus generally ACTRW and CTRW exhibit different behaviors.

More concretely, ACTRW describes the following process: a walker is trapped on the origin for time $t_1$, then jumps to $x_1$;  the walker is further trapped on $x_1$ for time $t_2$, and then jumps to a new position; this process is then renewed. Thus,  ACTRW process is characterized by a set of waiting times $\{t_1, t_2, \cdots, t_n, \cdots \}$ and displacements $\{x_1, x_2, \cdots, x_n, \cdots \}$. Here $\omega(t_a, t_1)$ is the PDF of the first waiting time $t_1$. In ACTRW process, the random walk starts from the time $t_a$, therefore $\omega(t_a, t_1)$ may depend on the aging time of the process $t_a$. The waiting times ${t_n}$ with $n>1$ are independent and identically distributed (i.i.d.) with a common probability density $\varphi(t)$. And the jump lengths $\{x_1, x_2, \cdots , x_n, \cdots \}$ are i.i.d. random variables, described by the probability density $f(x)$.

When $t_a=0$, we have $\omega(t_a, t_1)=\varphi(t_1)$, which is just the well known Montroll-Weiss nonequilibrium process. In order to investigate ACTRW, we should first discuss the aging renewal process. In what follows, we suppose that $P_{n_a}(t_a, t)$ is the probability of the renewal process $n_a(t_a, t)$, where $n_a(t_a, t)=n(t_a+t)-n(t)$ and $n(t)$ denotes the number of renewals by time $t$, i.e., $n_a(t_a, t)$ is the number of renewals in time interval $(t_{a},t_{a}+t)$ for a precess starts at the time $0$. Our main work is to discuss the properties of the renewal process \cite{Luck:1,Schulz:1} in the time interval $(t_a, t_a+t)$.

According to the renewal theory developed by Gord\`{e}che and Luck \cite{Luck:1},
\begin{equation} \label{ageq10}
\omega(s, u)=\frac{1}{1-\varphi(s)}\frac{\varphi(s)-\varphi(u)}{u-s},
\end{equation}
where $\omega(s, u)$ is the double Laplace transform %\cite{Dyke:1}
of the PDF of the first waiting time $\omega(t_a, t_1)$, and $\varphi(u)$ is the Laplace transform of $\varphi(t)$.
This paper focuses on taking $\varphi(t)$ as Eq. (\ref{ageq11}), and its Laplace transform  ($t\rightarrow u$) has the asymptotic form
\begin{equation}\label{ageq12}
\varphi(u)=e^{-(u+\lambda)^\alpha+\lambda^\alpha} \sim 1+\lambda^\alpha-(u+\lambda)^\alpha, \indent 0<\alpha<1;
\end{equation}
and if $\alpha=1$, $\varphi(u)\sim 1-u$.

The double Laplace transform of the PDF of  $n_a(t_a, t)$ reads \cite{Barkai:1}, $t_a\rightarrow s$, $t\rightarrow u$,
\begin{equation}\label{ageq13}
P_{n_a}(s, u)=\left\{ \begin{array}
 {l@{\quad } l}
 \frac{1-s\omega(s, u)}{su}, & n_a=0, \\
  \\ \omega(s, u)\varphi^{n_a-1}(u)\frac{1-\varphi(u)}{u} , &n_a\geq1. \\
 \end{array}
 \right.
\end{equation}
For the particular case $\alpha=1$, from (\ref{ageq12}),  $\omega(s, u)\sim 1/s$; and from Eq. (\ref{ageq13}) we can get $P_{n_a}(t_a,t)\sim \delta(n_a-t)$, for this case $P_{n_a}(t_a,t)$ is independent of $t_a$. Since $\omega(s, u)$ plays a key role in our discussion,  we now derive the analytical formulation of $\omega(t_a, t)$ \cite{Klafter:1},
\begin{equation}\label{ageq125}
\begin{split}
  \omega(s, t) & =\frac{1}{1-\varphi(s)} \exp(st)\int_t^{\infty}\exp(-sy)\varphi(y)dy \\
    & \sim \frac{(s+\lambda)^\alpha}{\lambda^\alpha-(s+\lambda)^\alpha}\frac{\exp(st)}{\Gamma(-\alpha)} \Gamma(-\alpha,(s+\lambda)t),
\end{split}
\end{equation}
%    &\sim [\frac{\lambda^\alpha}{\lambda^\alpha-(s+\lambda)^\alpha}-1]\frac{\exp(st)}{\Gamma(-\alpha)} \Gamma(-\alpha,(s+\lambda)t)
where $\Gamma(\alpha,x)=\int_x^{\infty}\exp(-t)t^{\alpha-1}d t$ is an incomplete Gamma function. Using the Laplace transform of incomplete Gamma function \cite{Oberhettinger:1}, we have,
\begin{equation}\label{ageq125a2}
\omega(t_{a}, t)=\frac{\exp(-\lambda t)}{-\Gamma(-\alpha)}g(t_{a})*(\exp(-\lambda*t_{a})(t_{a}+t)^{-\alpha-1}),
\end{equation}
where $g(t_{a})=t_{a}^{\alpha-1}\exp(-\lambda t_{a})E_{\alpha, \alpha}(\lambda^{\alpha}t_{a}^{\alpha})$.
From the second  line of Eq. (\ref{ageq125}), if $s\ll \lambda$, i.e., $t_{a}\gg 1/\lambda$. Eq. (\ref{ageq125}) can be given by
\begin{equation}\label{ageq125a1}
\omega(t_{a}, t) \sim \frac{\lambda \exp(-\lambda t)\sin(\pi \alpha)}{\alpha \pi t^\alpha}\int_{0}^{t_{a}}\exp(-\lambda \tau)\frac{1}{\tau+t}\tau^\alpha d \tau.
\end{equation}
%
%
%
%
%For $s \ll \lambda$, i.e. $t_{a} \gg 1/\lambda$, we can notice that $ \frac{\lambda^\alpha}{(s+\lambda)^\alpha-\lambda^\alpha} \gg 1$, Utilizing the Laplace transform property of Mittag-Leffer function $\mathcal{L}^{-1}\left\{ \frac{1}{(u+\lambda)^\alpha-\lambda^\alpha};t\right\}=t^{\alpha-1}e^{-\lambda t}E_{\alpha, \alpha}((\lambda t)^\alpha)$ leads to
%\begin{equation}\label{ageq126}
%\begin{split}
%  \omega(t_a, t) \sim & \frac{\lambda^\alpha \sin(\pi\alpha)\exp(-\lambda t)}{\pi}[t_a^{\alpha-1}\exp(-\lambda t_a)E_{\alpha,\alpha}(\lambda^\alpha t_a^\alpha)] \\
%  & *\left[\exp(-\lambda t_a)\frac{1}{t+t_a}\left(\frac{t_a}{t}\right)^\alpha\right],
%\end{split}
%\end{equation}
%where the  $*$ represents the convolution of functions.
If $s\gg \lambda$, i.e., $t_{a} \ll 1/\lambda$, then there exists $ \frac{(s+\lambda)^\alpha}{(s+\lambda)^\alpha-\lambda^\alpha} \rightarrow 1$. Eq. (\ref{ageq125}) can be further simplified as
%it can be noticed that the last line is the dominant term of Eq. (\ref{ageq126}), from Eq. (\ref{ageq125}) and (\ref{ageq126}) it can be noticed that Eq. (\ref{ageq126}) can be simplify to
\begin{equation}\label{ageq131}
  \omega(t_a, t) \sim \frac{\sin(\pi\alpha)\exp(-\lambda(t+t_a))}{\pi}\frac{1}{t+t_a}\left(\frac{t_a}{t}\right)^{\alpha};
\end{equation}
under the further assumption $t\ll 1/\lambda$, i.e.,  $\lambda t\ll 1$, there exists
%\textcolor[rgb]{1.00,0.00,0.00}{when both $t\ll 1/\lambda$ and $t_{a}\ll 1/\lambda$ , i.e.,  $\lambda t\ll 1$ and $\lambda t_{a}\ll 1$ we can get
\begin{equation}\label{ageq132}
\begin{split}
  \omega(t_a, t) \sim & \frac{\sin(\pi\alpha)}{\pi}\frac{1}{t+t_a}\left(\frac{t_a}{t}\right)^{\alpha},
\end{split}
\end{equation}
being the same as the one given in \cite{Klafter:1,Barkai:1} for the power law waiting time, i.e., $\lambda=0$.
%this result is consistent with when the waiting time $\varphi(t)$ is power law, i.e., $\lambda=0$, the first term of Eq. (\ref{ageq126}) is omitted, being same as the one given in \cite{Klafter:1,Barkai:1}.
 When $\lambda$ is sufficiently large, Eq. (\ref{ageq125a1}) plays a dominant role. %}
 \begin{figure}[htb]
  \centering
  \includegraphics[width=9cm, height=6cm]{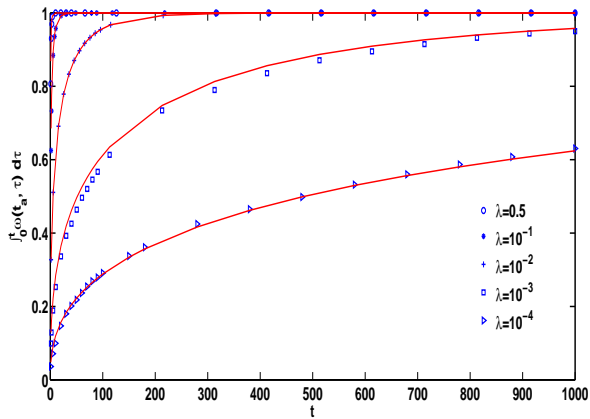}\\
  \caption{Probability of particles making jumps during the time interval $[t_a,t_a+t]$ for $t_{a} \gg 1/\lambda$.
  The parameters of $\varphi(t)$ Eq. (\ref{ageq11}) are taken as $\alpha=0.6$, $t_{a}=10^{4}$, and $t=1000$; and the symboled lines are obtained by averaging $5000$ trajectories with different $\lambda$. The solid lines from down to up corresponding to the increased $\lambda$ are the theoretical results of Eq. (\ref{ageq125a2}). When $\lambda$ is large, the probability reaches $1$ quickly.  }\label{agfig15}
\end{figure}
\begin{figure}[htb]
  \centering
  \includegraphics[width=9cm, height=6cm]{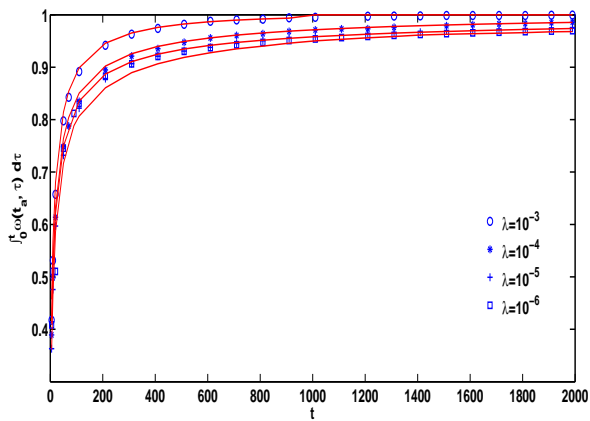}\\
  \caption{Probability of particles making jumps during the time interval $[t_a,t_a+t]$.
  The parameters  are taken as $\alpha=0.6$, $t_{a}=100$, and $t=2000$; and the symboled lines are obtained by averaging $5000$ trajectories with different $\lambda$. The solid lines from down to up corresponding to the increased $\lambda$ are the theoretical results of Eq. (\ref{ageq125a2}). When $\lambda$ is sufficiently small, the distribution is almost the same as pure power law for short times.  }\label{agfig25}
\end{figure}
In the following, we analyze the asymptotic form of $\omega(t_a, t)$ with $\alpha \in (0,1)$.
From  Eq. (\ref{ageq10}) and Eq. (\ref{ageq12}), there exists
\begin{equation}\label{ageq14}
\omega(s, u)\sim \frac{(s+\lambda)^\alpha-(u+\lambda)^\alpha}{(u-s)[\lambda^\alpha-(s+\lambda)^\alpha]}.
\end{equation}
For Eq. (\ref{ageq14}), taking the limit $s\rightarrow \infty$, we have $\omega(s, u) \sim s^{-1}$; then $\omega(0, t) \sim \delta(t)$. Using the relation $P_{0}(t_{a}, t)=1-\int_{0}^{t}\omega(t_{a}, t) d t$ leads to $P_{0}(0, t)=0$, i.e., all the particles move at time $t=0$, no aging phenomenon.

\begin{figure}[htb]
  \centering
  \includegraphics[width=9cm, height=6cm]{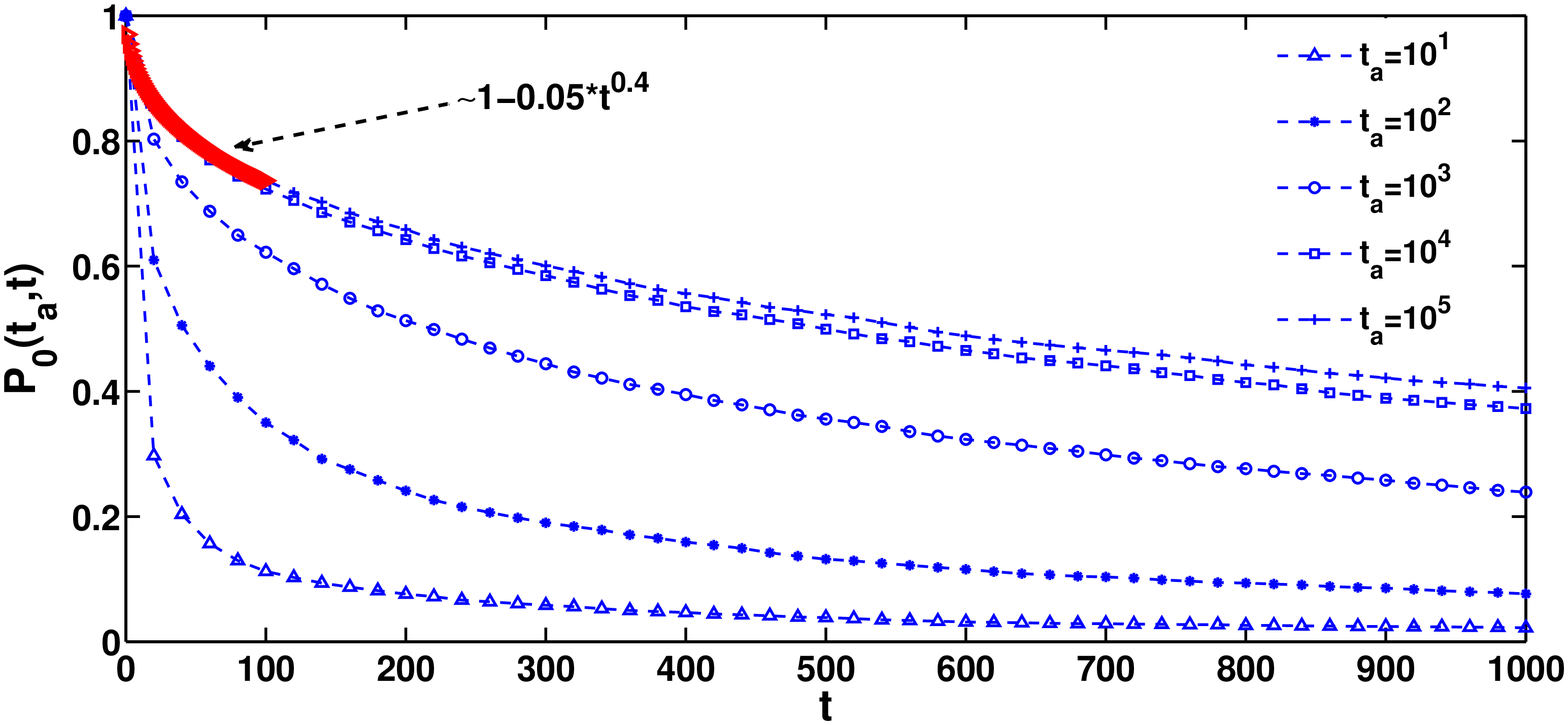}\\
  \caption{Time evolution of $P_{0}(t_{a},t)$ with different $t_a$. The parameters are taken as $\alpha=0.6$, $\lambda=10^{-4}$, and $t=1000$.  The lines are obtained by averaging $5000$ trajectories. The (red) dashed line, $1-0.05*t^{0.4}$, is the fitting result for small $t$ and big $t_a$, which agrees with Eq. (\ref{ageq133}).
%  \textcolor[rgb]{1.00,0.00,0.00}{the red  dash line is the fitting result for small $t$ when $t_{a}$ is large, our fitting result is $1-0.05*t^{0.4}$, it's is agree with Eq. (\ref{ageq133}) }
  }\label{agfig20}
\end{figure}

%\textcolor[rgb]{1.00,0.00,0.00}{From Eq. (\ref{ageq14}), it can be noticed that $\omega(t_{a}, t)$ is independent of $t_{a}$ in the limit $t_{a}\rightarrow 0$. Consider the limit $s\rightarrow \infty$, we can find $\omega(s, u) \sim s^{-1}$, then we have $\omega(0, t) \sim \delta(t)$, using the relation  $P_{0}(t_{a}, t)=1-\int_{0}^{t}\omega(t_{a}, t) d t$, we have $P_{0}(0, t)=0$, i.e., all the particle move at time $t=0$, there is no aging phenomenon,}
Consider the survival probability $P_{0}(t_a, t)$ \cite{krusemann:3}, which gives the probability of making no jumps during the interval $t_a$ up to $t_a+t$,
\begin{equation}\label{ageq15}
P_0(s, u)\sim\frac{1}{su}-\frac{(s+\lambda)^\alpha-(u+\lambda)^\alpha}{u(u-s)[\lambda^\alpha-(s+\lambda)^\alpha]}.
\end{equation}
It is instructive to consider two different limits. If $s\ll u$, i.e., $t_a\gg t$, there exists
\begin{equation}\label{ageq130}
  P_0(s, u)\sim \frac{1}{su}-\frac{(u+\lambda)^\alpha-\lambda^\alpha}{u^2[(s+\lambda)^\alpha-\lambda^\alpha]}.
\end{equation}
 For $\lambda \ll u$, i.e., $ t_{0} \ll t\ll 1/\lambda$, then $(u+\lambda)^\alpha \sim u^\alpha(1+\lambda/u)^\alpha \sim u^\alpha$. Performing double inverse Laplace transform on the above equation results in
\begin{equation}\label{ageq16}
\begin{split}
  P_0(t_a, t)  \sim & \, 1-t_a^{\alpha-1}\exp(-\lambda t_a)E_{\alpha, \alpha}(\lambda^\alpha t_a^\alpha)\Big(-\lambda^\alpha t
  \\
    & +\frac{t^{1-\alpha}}{\Gamma(2-\alpha)} \Big).
\end{split}
\end{equation}
 For $t_a\ll 1/\lambda$, Eq. (\ref{ageq16}) can be simplified as
\begin{equation}\label{ageq124}
  P_0(t_a, t)  \sim  1- t_a^{\alpha-1}\Big(-\lambda^\alpha t+\frac{t^{1-\alpha}}{\Gamma(2-\alpha)}\Big).
\end{equation}
It can be noted that $t^{1-\alpha}/\Gamma(2-\alpha)$ is larger than $\lambda^\alpha t$ in the parenthesis of Eq. (\ref{ageq124}), since $\lambda^\alpha t =(\lambda t)^\alpha t^{1-\alpha} \ll t^{1-\alpha}$. For $t_a\gg 1/\lambda$, from Eq. (\ref{ageq16}) and Eq. (\ref{age918}), we obtain
\begin{equation}\label{ageq133}
  P_0(t_a, t)  \sim  1-\frac{t^{1-\alpha}}{\langle\tau\rangle\Gamma(2-\alpha)},
\end{equation}
being confirmed by FIG. \ref{agfig20}, i.e., the lines tend to be close for big $t_a$.

%\textcolor[rgb]{1.00,0.00,0.00}{while for $t_a\gg 1/\lambda$, from Eq. (\ref{ageq16}) and Eq. (\ref{age918}) we obtain
%\begin{equation}\label{ageq133}
%  P_0(t_a, t)  \sim  1-\frac{t^{1-\alpha}}{\langle\tau\rangle\Gamma(2-\alpha)}
%\end{equation}
%From the Eq. (\ref{ageq124}), it can noticed that with $P_0(t_a, t)$ is independent of $t_{a}$ as $t_{a}$ grow. the results is agree with FIG. \ref{agfig20}.}

And for $t\ll t_a\ll 1/\lambda$, we have $P_0(t_{a}, t) \sim 1$, i.e., for small $\lambda$, $\varphi(t) \sim t^{-1-\alpha}$, the waiting time is generally long; in a small observation time $t$, we cannot find movement of the particles.
Eq. (\ref{ageq15}) can also be rewritten as
\begin{equation}\label{ageq127}
 P_{0}(s,u) \sim \frac{1}{su}+\frac{1}{u(u-s)}+\frac{(u+\lambda)^\alpha-\lambda^\alpha}{u(s-u)[(s+\lambda)^\alpha-\lambda^\alpha]}. \end{equation}
For the case $u \ll s $, i.e., $t_{a}\ll t$,  Eq. (\ref{ageq127}) yields
\begin{equation}\label{ageq128}
 P_{0}(s,u) \sim \frac{(u+\lambda)^\alpha-\lambda^\alpha}{us[(s+\lambda)^\alpha-\lambda^\alpha]};
\end{equation}
under the further assumption $u \gg \lambda$, we have
\begin{equation}\label{ageq129}
 P_{0}(t_{a},t)\sim 1*g(t_{a})\frac{t^{-\alpha}}{\Gamma(1-\alpha)},
\end{equation}
where $g(t_{a})$ is defined in Eq. (\ref{age917}), i.e., when $t_a \ll t \ll 1/\lambda$,
\begin{equation}\label{ageq1291}
 P_{0}(t_{a},t) \sim \frac{\sin(\pi\alpha)}{\pi\alpha}\left(\frac{t}{t_{a}}\right)^{-\alpha},
\end{equation}
being the same as the result given in \cite{Luck:1} for the pure power law case ($\lambda=0$).

%\textcolor[rgb]{1.00,0.00,0.00}{while for $t \gg 1/\lambda$, from Eq. (\ref{ageq130}) we have $p_{0}(t_{a},t)=\frac{\langle \tau^2\rangle}{\langle \tau\rangle} \delta(t)$.}

%\textcolor[rgb]{1.00,0.00,0.00}{In the following we consider the case for $t_{a}\ll t$,  Eq. (\ref{ageq15}) can be given in another way,
%%\begin{equation}\label{ageq127}
%% P_{0}(s,u)=\frac{1}{su}-\frac{1}{u(u-s)}+\frac{(u+\lambda)^\alpha-\lambda^\alpha}{u(s-u)(s+\lambda)^\alpha-\lambda^\alpha},
%%\end{equation}
%for $u\ll s$, Eq. (\ref{ageq127}) yields,
%\begin{equation}\label{ageq128}
% P_{0}(s,u)=\frac{(u+\lambda)^\alpha-\lambda^\alpha}{us(s+\lambda)^\alpha-\lambda^\alpha},
%\end{equation}
%for $u \gg \lambda$, we have
%\begin{equation}\label{ageq129}
% P_{0}(t_{a},t)\sim 1*g(t_{a})\frac{t^{-\alpha}}{\Gamma(1-\alpha)},
%\end{equation}}
%
%
%\textcolor[rgb]{1.00,0.00,0.00}{as $t\ll 1/\lambda$, we have $t_{a} \ll 1/\lambda$,
%\begin{equation}\label{ageq1291}
% P_{0}(t_{a},t) \sim \frac{\sin(\pi\alpha)}{\pi\alpha}\left(\frac{t}{t_{a}}\right)^{-\alpha}.
%\end{equation}
%the result has been shown in \cite{Luck:1}.}
\begin{figure}[htb]
  \centering
  \includegraphics[width=9cm, height=6cm]{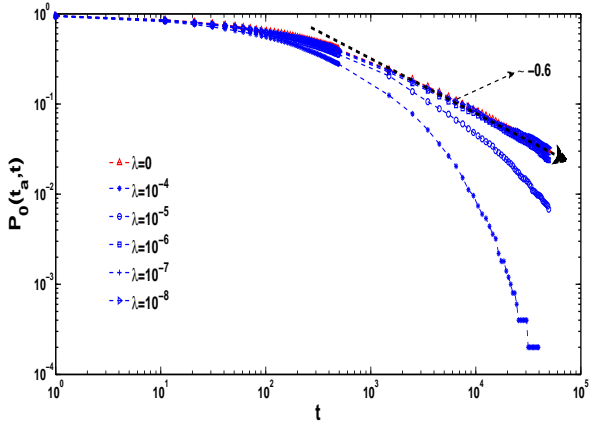}\\
  \caption{Time evolution of $P_{0}(t_{a},t)$ with different $\lambda$.
  %the waiting time PDF Eq. (\ref{ageq11}) for $t \ll 1/\lambda$  and $t_{a} \ll 1/\lambda$.
  The parameters are taken as $\alpha=0.6$, $t_{a}=500$, and $t=50000$.  The lines are obtained by averaging $5000$ trajectories. The dashed line with arrow is for the indicator of slope $-0.6$, confirming Eq. (\ref{ageq1291}).
  %When $t \gg t_{a}$ we can get $P_{0}(t_{a},t) \sim \sin(\pi \alpha)/\pi/\alpha (t/t_{a})^{-\alpha}$ \cite{Luck:1} with waiting time is power law, our fitting result is consistent with $P_{0}(t_{a},t) \sim t^{-\alpha}$.
  }\label{agfig24}
\end{figure}
When the aging time is sufficiently long relative to the observation time and $\lambda$ is small, the probability of making no jumps during the time interval $(t_a, t_a+t)$ approaches to one, i.e., the system is completely trapped. On the contrary, if $t_a$ is short, while the observation time is long enough, then the particles are unacted on the aging time. So, at least one jump will be made, namely, the possibility of making no jumps is zero. Indeed, from Eq. (\ref{ageq1291}), it can be easily obtained that $P_0(t_{a}, t)\sim 0$ when $t/t_{a}\rightarrow \infty$; and $P_0(t_{a}, t)\sim 0$ can also be directly obtained from Eq. (\ref{ageq15}) under the assumption $t \gg t_a$.

%From Eq. (\ref{ageq1291}), we have when $t/t_{a}\rightarrow \infty$ we have $P_0(t_{a}, t)\sim 0$, the result can be get from Eq. (\ref{ageq15}) directly.

%These results are apparent, when the aging time is enough long, the probability of making no steps during the time interval $(t_a, t_a+t)$ is approaches to one, namely, it means no steps happen when $t_a$ is long enough, because the system becomes completely trapped, the particles are almost motionless. On the contrary, when $t_a$ is a short time while our observation time is long  enough relatively, then the particles are unacted on the aging time, so at least one step will happen at least, therefore, we have $P_0(t_a, t)\sim 0$, namely, the possibility of making no steps is zero. \\
From Eq. (\ref{ageq13}), we can write the double Laplace transform of the PDF of $n_a(t_a, t)$ as
\begin{equation}\label{ageq17}
P_{n_a}(s, u)=\frac{\delta(n_a)}{u}\left[\frac{1}{s}-\omega(s, u)\right]+\omega(s, u)\varphi^{n_a-1}(u)\frac{1-\varphi(u)}{u}.
\end{equation}
Inserting Eq. (\ref{ageq12}) into the above equation yields
\begin{equation}\label{ageq18}
\begin{split}
P_{n_a}(s, u)  \sim  &  ~    \omega(s, u)\exp^{-n_a[(u+\lambda)^\alpha-\lambda^\alpha]}\frac{(u+\lambda)^\alpha-\lambda^\alpha}{u} \\
    & +\frac{\delta(n_a)}{u}\left[\frac{1}{s}-\omega(s, u)\right].
\end{split}
\end{equation}
From now on, we start to calculate the $p$-th moment of the aging renewal process $n_a(t_a, t)$ \cite{note:1}, which reads
\begin{equation}\label{ageq19}
\langle n_a^p(s, u)\rangle=\frac{\Gamma(p+1)\omega(s, u)}{ u[(u+\lambda)^\alpha-\lambda^\alpha]^p}.
\end{equation}
For the cases that $u\gg s$ and $\lambda \ll u$ or $u \ll s$, there exist
\begin{equation}\label{ageq115}
\langle n_a^p(s, u)\rangle \sim \left\{
                            \begin{array}{ll}
                             \frac{\Gamma(p+1)}{((s+\lambda)^{\alpha}-\lambda^{\alpha})u^{2+\alpha(p-1)}},  & \hbox{for } s\ll u,\\
                             & ~~ ~~ \lambda \ll u; \\
                              \frac{\Gamma(p+1)}{(\alpha \lambda^{\alpha-1})^{p}}\frac{1}{s u^{1+p}},  & \hbox{for } u\ll s.
                            \end{array}
                          \right.
\end{equation}
By double inverse Laplace transform we have
\begin{equation}\label{ageq116}
\langle n_a^p(t_{a}, t)\rangle \sim \left\{
                            \begin{array}{ll}
                              \frac{\Gamma(p+1)}{\Gamma(2+\alpha p-\alpha)}g(t_{a})t^{\alpha p-\alpha +1},  & \hbox{for } t\ll t_{a}, \\
                              & ~~ ~~ \lambda t \ll 1; \\
                              (\frac{t}{\langle \tau \rangle})^{p},  & \hbox{for } t_{a}\ll t.
                            \end{array}
                          \right.
\end{equation}
 Taking $p=1$ in (\ref{ageq19}) leads to
\begin{equation}\label{ageq110}
\langle n_a(s, u)\rangle=\frac{(s+\lambda)^\alpha-(u+\lambda)^\alpha}{u(u-s)[\lambda^\alpha-(s+\lambda)^\alpha][(u+\lambda)^\alpha-\lambda^\alpha]},
\end{equation}
which can be rewritten as
\begin{equation}\label{ageq111}
\langle n_a(s, u)\rangle =\frac{1}{u(s-u)}\left(\frac{1}{(u+\lambda)^\alpha-\lambda^\alpha}-\frac{1}{(s+\lambda)^\alpha-\lambda^\alpha}\right).
\end{equation}
We will confirm that if both $t_{a}$ and $t$ are large scales, $ \langle n_a(t_a, t)\rangle\sim t/\langle\tau\rangle$, which is an important result for normal diffusion. For small $u$ and $s$, using the Taylor expansion $(u+\lambda)^{\alpha} \sim \lambda^\alpha+\alpha\lambda^{\alpha-1}u$ and $(s+\lambda)^{\alpha} \sim \lambda^\alpha+\alpha\lambda^{\alpha-1}s $,  from Eq. (\ref{ageq111}), we have
\begin{equation}\label{ageq120}
\langle n_a(s, u)\rangle \sim  \frac{1}{\alpha \lambda^{\alpha-1}u^2s}.
\end{equation}
Performing double Laplace transform on the above equation yields
\begin{equation}\label{ageq121}
\langle n_a(t_{a}, t)\rangle \sim  \frac{t}{\alpha\lambda^{\alpha-1}}=\frac{t}{\langle \tau\rangle},
\end{equation}
where $\langle \tau\rangle=\alpha\lambda^{\alpha-1}$ (see Appendix \ref{AppendixB}).

For the slightly aging system, $t\gg t_a$, i.e., $u\ll s$; performing the double inverse Laplace transform on both sides of (\ref{ageq111}) yields
\begin{equation}\label{ageq112}
\begin{split}
 \langle n_a(t_a, t)\rangle \sim & (t^{\alpha-1}e^{-\lambda t}E_{\alpha, \alpha}(\lambda^\alpha t^{\alpha}))*1 \\
 & -(t_{a}^{\alpha-1}e^{-\lambda t_{a}}E_{\alpha, \alpha}(\lambda^\alpha t_{a}^{\alpha}))*1 \\
     \sim & (t^{\alpha-1}e^{-\lambda t}E_{\alpha, \alpha}(\lambda^\alpha t^{\alpha}))*1.
\end{split}
\end{equation}
For the special case, $\lambda=0$,  it can be noted that $\langle n_a(t_a, t)\rangle \sim t^{\alpha}$.  %In our simulation we suppose the relation $\langle n_a(t_a, t)\rangle \sim t^{\beta} t_{a}^{\gamma}$.
For $t\gg 1/\lambda$, using the asymptotic expansion of Mittag-Leffler function (\ref{ageq915}), from Eq. (\ref{ageq112}), we again obtain
 $\langle n_a(t_a, t) \rangle \sim \frac{t}{\langle\tau\rangle}$.
%where $ \langle \tau \rangle =\Gamma(-\alpha)\lambda^{\alpha-1}$,
In the long time scale, the process converges to the Gaussian process, and then the first moment of the number of renewal events grows linearly with the observation time $t$. For $t \ll 1/\lambda$, from Eq. (\ref{ageq112}) we have
%the asymptotic expansion of the Mittag-Leffler function (\ref{age913}) leads to \cite{Luck:1}
\begin{equation}\label{ageq123}
 \langle n_a(t_a, t) \rangle \sim \frac{1}{\Gamma(1+\alpha)}t^{\alpha}.
\end{equation}
It can be seen that when $t\gg t_a$ the first moment of $n_a$ is not relevant to the aging time $t_a$. From Eq. (\ref{ageq121}) and  Eq. (\ref{ageq123}), we can see that $\lambda$ plays an important role in our discussion as expected.

\begin{figure}[htb]
  \centering
  \includegraphics[width=9cm, height=6cm]{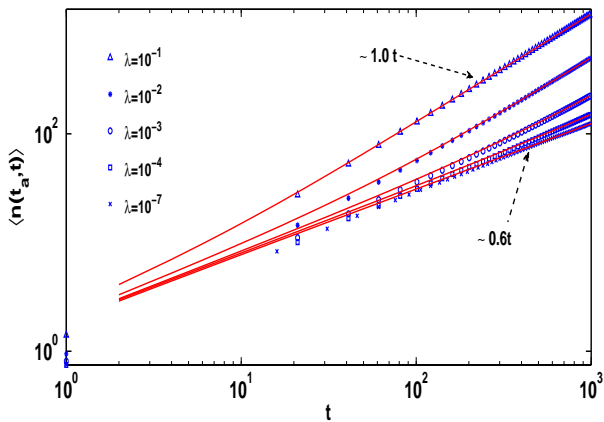}\\
  \caption{Time evolution of the ensemble average of the renewal times $\langle n_a(t_a,t)\rangle$ with the waiting time PDF Eq. (\ref{ageq11}) for slight aging. The parameters are taken as $\alpha=0.6$, $t_a=3$, and $t=1000$. The real lines are for the analytical result Eq. (\ref{ageq112}) and the other lines are obtained by averaging $5000$ trajectories.
%   It can been that when time $t$ is short, the lines are straight. And $\lambda$ is smaller, the straight line is longer.
 % Simulation of $\langle n_a(t_a,t)\rangle$ for the PDF of Eq. (\ref{ageq11}) for slightly aging. The index of $\beta$ for different $\lambda$,   drawn for $\lambda=0. 1,
%0. 05, 0. 03, 0. 02$ and $10^{-10}$.  the parameter $\beta$ ranges from $0. 6$ to $1$. Normal diffusion for $\beta=1$(larger$\lambda$),   and subdiffusion for $\beta=0. 6$ (smaller $\lambda$).  The parameters are $\alpha=0. 6$,  $t_a=3$,  $t=500$.
}\label{agfig12}
\end{figure}

%\begin{figure}[htb]
%  \centering
%  % Requires \usepackage{graphicx}
%  \includegraphics[width=9cm, height=6cm]{agfig91.eps}\\
%  \caption{Time evolution of the ensemble average of the renewal times $\langle n_a(t_a,t)\rangle$ with the waiting time PDF Eq. (\ref{ageq11}) for slight aging. The parameters are taken as $\alpha=0.6$, $t_a=3$, and $t=1000$. The real lines are for the analytical result Eq. (\ref{ageq112}) and the other lines are obtained by averaging $5000$ trajectories. It can been that when time $t$ is short, the lines are straight. And $\lambda$ is smaller, the straight line is longer.
% % Simulation of $\langle n_a(t_a,t)\rangle$ for the PDF of Eq. (\ref{ageq11}) for slightly aging. The index of $\beta$ for different $\lambda$,   drawn for $\lambda=0. 1,
%%0. 05, 0. 03, 0. 02$ and $10^{-10}$.  the parameter $\beta$ ranges from $0. 6$ to $1$. Normal diffusion for $\beta=1$(larger$\lambda$),   and subdiffusion for $\beta=0. 6$ (smaller $\lambda$).  The parameters are $\alpha=0. 6$,  $t_a=3$,  $t=500$.
%}\label{agfig12}
%\end{figure}

While for the strongly aging system, $t_a\gg t$, i.e., $s \ll u$; there exists
\begin{equation}\label{ageq113}
\langle n_a(s, u)\rangle \sim\frac{1}{u^2[(s+\lambda)^\alpha-\lambda^\alpha]},
\end{equation}
which yields
\begin{equation}\label{ageq114}
\langle n_a(t_a, t)\rangle \sim tt_a^{\alpha-1}e^{-\lambda t_a}E_{\alpha, \alpha}(\lambda^\alpha t_a^\alpha).
\end{equation}
%It shows that the aging time $t_a$ plays an important role.
Following the way used above, for $t_{a}\gg 1/\lambda$, the term $t_a^{\alpha-1}E_{\alpha, \alpha}(\lambda^\alpha t_a^\alpha)e^{-\lambda t_a}$ tends to $\lambda^{1-\alpha}/\alpha$ (\ref{age918}). Then we have
\begin{equation}\label{ageq117}
\langle n_a(t_a, t)\rangle \sim \frac{t}{\langle\tau\rangle}.
\end{equation}
%The term $t_a^{\alpha-1}E_{\alpha, \alpha}(\lambda^\alpha t_a^\alpha)e^{-\lambda t_a}$ tends to $\lambda^{1-\alpha}/\alpha$, for large scale $t_a$ lose its role in the $\langle n_a(t_a, t)\rangle$.
For $t_{a} \ll 1/\lambda$, there exists \cite{Luck:1,Klafter:1}
\begin{equation}\label{ageq118}
\langle n_a(t_a, t)\rangle \sim \frac{1}{\Gamma(\alpha)}t t_{a}^{\alpha-1}.
\end{equation}

The above results for the first moment of $n_a(t_a, t)$ can be summarized as: 1. if $t_a$ or $t$ is greater than $1/\lambda$, then $\langle n_a(t_a, t)\rangle \sim t/\langle \tau \rangle$; 2. for $t\gg t_a$ and $t\ll1/\lambda$ (i.e. $t^{-\alpha-1}\exp(-\lambda t)\sim t^{-\alpha-1}$), $\langle n_a(t_a, t)\rangle \sim t^\alpha/\alpha$; 3. for $t\ll t_{a}$ and  $t_a \ll 1/\lambda$, $\langle n_a(t_a,t)\rangle$ behaves as $t t_{a}^{\alpha-1}$.

%For the case of $0<\alpha<1$, we can find that,
%1, each of $t_a$ and $t$ great than $1/\lambda$, $\langle n_a(t_a, t)\rangle \sim t/\langle \tau \rangle$
%2, for $t\gg t_a$,and $t\ll1/\lambda$, $\langle n_a(t_a, t)\rangle \sim t^\alpha/\alpha$, i.e. $t^{-\alpha-1}\exp(-\lambda t) \sim t^{-\alpha-1}$, the result is similar to power law distribution.
%3, if $t\ll t_{a}$ and for $t_a \ll 1/\lambda$, $\langle n_a(t_a,t)\rangle$ has the behaviour $t t_{a}^{\alpha-1}$.

\begin{figure}[htb]
  \centering
  \includegraphics[width=9cm, height=6cm]{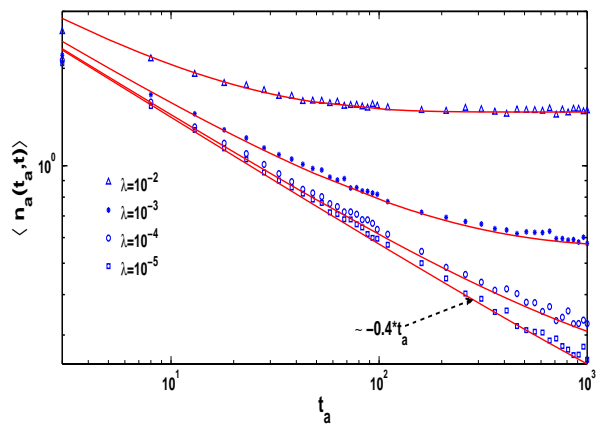}\\
  \caption{ The relation between the $\langle n_{a}(t_{a}, t)\rangle $ and $\lambda$ for $t_{a}\gg t$ according to the trajectories and the theory. The number of particles is $3000$,  $t_{a}=1000$,  $t=5$ , $\alpha=0. 6$,  and $\lambda=10^{-2}$, $10^{-3}$, $10^{-4}$, and $10^{-5}$. The real lines are for the analytical result Eq. (\ref{ageq114}) and the other symbol lines are obtained by averaging $5000$ trajectories.}\label{agfig13}
\end{figure}
%From Eq.(\ref{ageq117}), we know that $\gamma=-\alpha-3$, the dash-dot line with arrow is our fitting result, the formula is $\gamma=-0.3389 \log(\lambda)-2.306$. For $t_a\ll 1/lambda$, we have
\begin{figure}[htb]
  \centering
  \includegraphics[width=9cm, height=6cm]{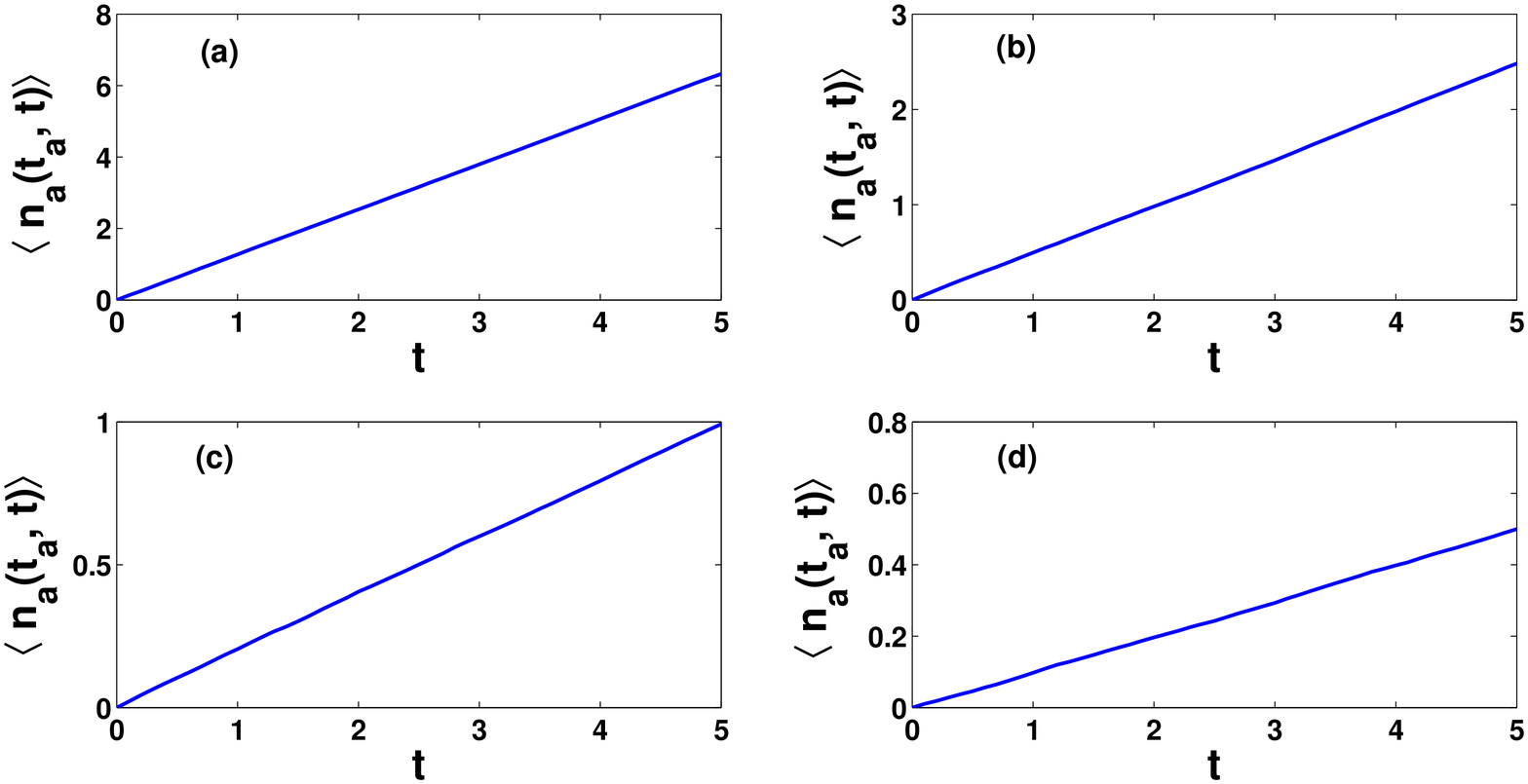}\\
  \caption{Time evolution of the ensemble average of the renewal times $\langle n_{a}(t_{a}, t)\rangle$ for strong aging ($t_{a} \gg t$). The parameters $\alpha=0.6$, $t_a=1000$, and $\lambda=10^{-1}$ for (a); $\lambda=10^{-2}$ for (b); $\lambda=10^{-3}$ for (c); $\lambda=10^{-4}$ for (d). The lines are obtained by averaging $10^4$ trajectories.  It can be seen that $\langle n_{a}(t_{a}, t)\rangle$ grows linearly with time $t$ for fixed $t_{a}$ for all  values of $\lambda$, which confirms the analytical result Eq. (\ref{ageq114}).
 % The relation between the $\langle n_{a}(t_{a}, t)\rangle  $ and  the observation time $t$ for $t_{a} \gg t$ getting from  the trajectories  of the particles.  From the figure we can find that $\langle n_{a}(t_{a}, t)\rangle$ grows linearly with $t$ for fixed $t_{a}$.
  }\label{agfig14}
\end{figure}

\section{Tempered ACTRW}\label{sec3}

\subsection{Mean squared displacement} \label{subsec3-1}

After understanding the statistics of the number of renewals, we go further to discuss the tempered ACTRW with the waiting time distribution Eq. (\ref{ageq11}). The process of ACTRW has been described in Sec. \ref{sec2}. This paper focuses on the symmetric random walk, i.e., the distribution of jump lengths $f(x)=f(-x)$; and $M_2=\int_{-\infty}^{+\infty} x^2 f(x) dx$ is finite. For such a random walk, we denote $P(x, t_a, t)$ as the PDF of particles' position in the decoupled tempered ACTRW with aging time $t_a$. Then
\begin{equation}\label{ageq301}
P(x, t_a, t)=\sum^{\infty}_{n_a=0}P_{n_a}(t_a, t)f_{n_a}(x),
\end{equation}
where $P_{n_a}(t_a, t)$ means the probability of jumping $n_a$ steps in the time interval $(t_a, t_a+t)$, and $f_{n_a}(x)$ the probability of jumping to the position $x$ after $n_a$ steps. In the Fourier-Laplace domain,
\begin{equation}\label{ageq302}
P(k, s, u)=\sum^{\infty}_{n_a=0}P_{n_a}(s, u)f^{n_a}(k).
\end{equation}
 Inserting Eq. (\ref{ageq13}) into Eq. (\ref{ageq302}) leads to
\begin{equation}\label{ageq303}
P(k, s, u)=\frac{1-s\omega(s, u)}{su}+\frac{\omega(s, u)(1-\varphi(u))}{u}\frac{f(k)}{1-\varphi(u)f(k)}.
\end{equation}
By differentiating Eq. (\ref{ageq303}) two times with respect to $k$ and setting $k=0$, we derive the second order moment of the random walks, i.e.,
\begin{equation}\label{ageq304}
\langle r^{2}(s, u)\rangle=\frac{\omega(s, u)M_2}{u[1-\varphi(u)]}.
\end{equation}
For the mean square displacement, we present the results of the slightly aging and strongly aging system, i.e.,
\begin{equation}\label{ageq305}
\langle r^{2}(s, u)\rangle \sim \left\{ \begin{array}
 {l@{\quad } l}
 \frac{M_2}{us[(u+\lambda)^\alpha-\lambda^\alpha]}, & \hbox{for } u\ll s, \\
  \\ \frac{M_2}{u^2[(s+\lambda)^\alpha-\lambda^\alpha]} , &\hbox{for } u\gg s. \\
 \end{array}
 \right.
\end{equation}
Performing double inverse Laplace transform on $\langle r^{2}(s, u)\rangle$ yields
\begin{equation}\label{ageq306}
\langle r^{2}(t_a, t)\rangle \sim \left\{ \begin{array}
 {l@{\quad } l}
M_{2}g(t)*1,  & \hbox{for } t\gg t_a,  \\
  \\ M_{2}t g(t_a), &  \hbox{for } t_a\gg t,  \\
 \end{array}
 \right.
\end{equation}
where $g(z)=z^{\alpha-1}\exp(-\lambda z)E_{\alpha, \alpha}(\lambda^{\alpha}z^{\alpha})$.

From FIG. (\ref{agfig22}), we can see the large fluctuations of $\langle r^{2}(t_a, t)\rangle$ even if the number of trajectories is 10000. This is because that most of particles are trapped in the initial position for $t_{a}\gg t \ll 1/\lambda$, which is consistent with Eq. (\ref{ageq16}). This is related to population splitting \cite{Cherstvya:1}.

It can be noted that when $t \gg t_a$, the mean squared displacement has no aging effect; while $t \ll t_a$ ($\lambda t_a \ll 1 $), it is deeply affected by the aging time $t_a$. The surprising result is that $\langle r^{2}(t_a,  t)\rangle \sim \langle n_a(t_a,  t)\rangle$, when the second order moment of the jump length is finite; the same things happen for the pure power-law waiting time distribution.

%This result is valid,  in the limit $t \gg t_a$,  we recover the system is nothing with the aging. In the aging regime,  $t\ll t_a$,  we can see a fact is that the system is deeply affected by the aging time $t_a$. and moreover especially,  we can find a surprising result,  that is the order of $\langle r^{2}(t_a,  t)\rangle$ is almost similar to the order of $\langle n_a(t_a,  t)\rangle$,  no matter which one the limit is,  this result all set up. This is also a same result for the pure power-law waiting time distribution situation.

\begin{figure}[htb]
  \centering
  \includegraphics[width=9cm,  height=6cm]{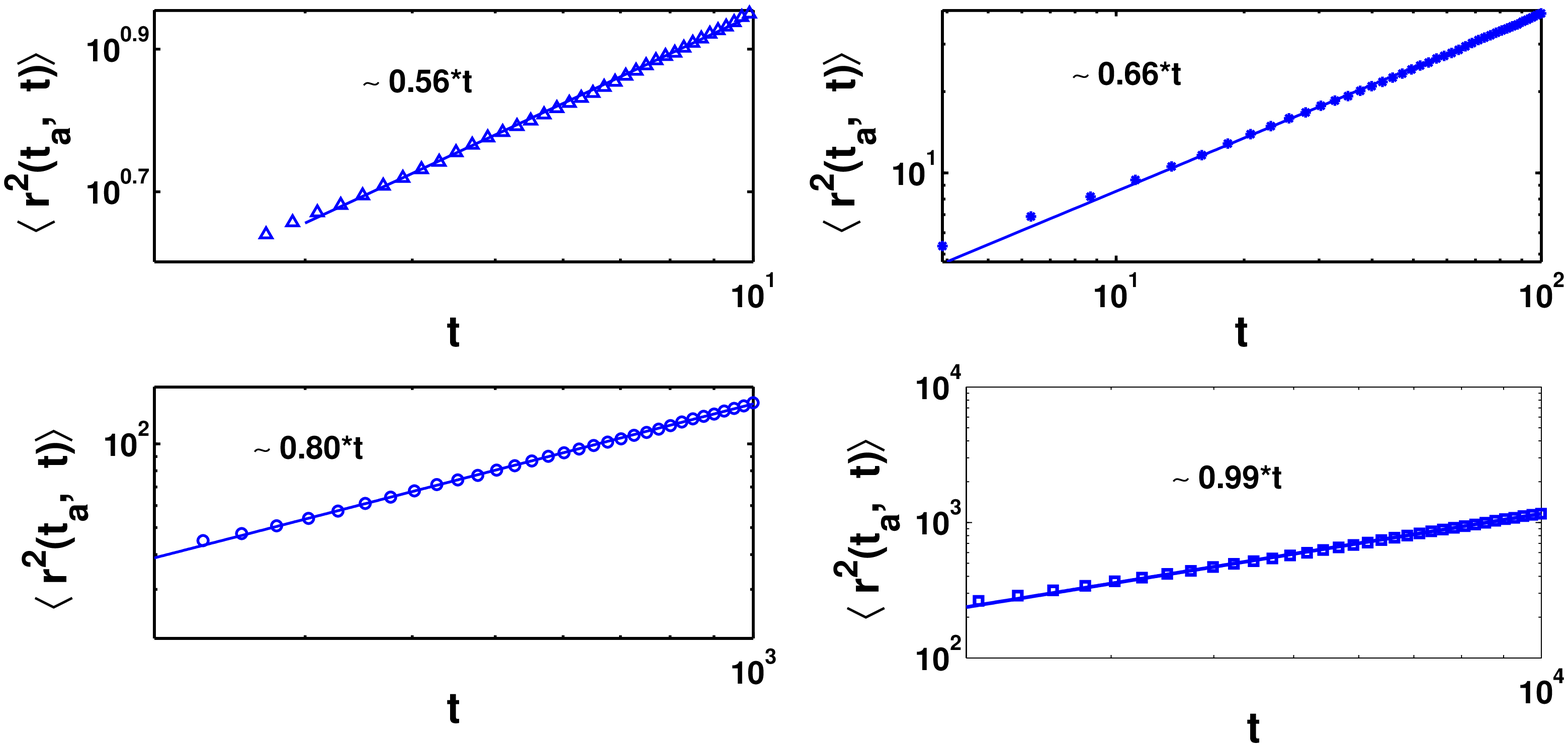}\\
  \caption{ The relation between the $\langle r^2(t_{a},  t)\rangle $ and the observation time $t$ for $t_{a} \ll t$ getting from  the trajectories of the particles (dashed line) and Eq. (\ref{ageq306}) (real line). The parameters $\alpha=0.6$, $\lambda=10^{-3},$ $t_a=1$, and the number of trajectories is $5000$. With the increase of $\lambda$, the characteristic of $\langle r^2(t_{a},  t)\rangle $ changes from Power law to normal diffusion.    }\label{agfig21}
\end{figure}

\begin{figure}[htb]
  \centering
  \includegraphics[width=9cm,  height=6cm]{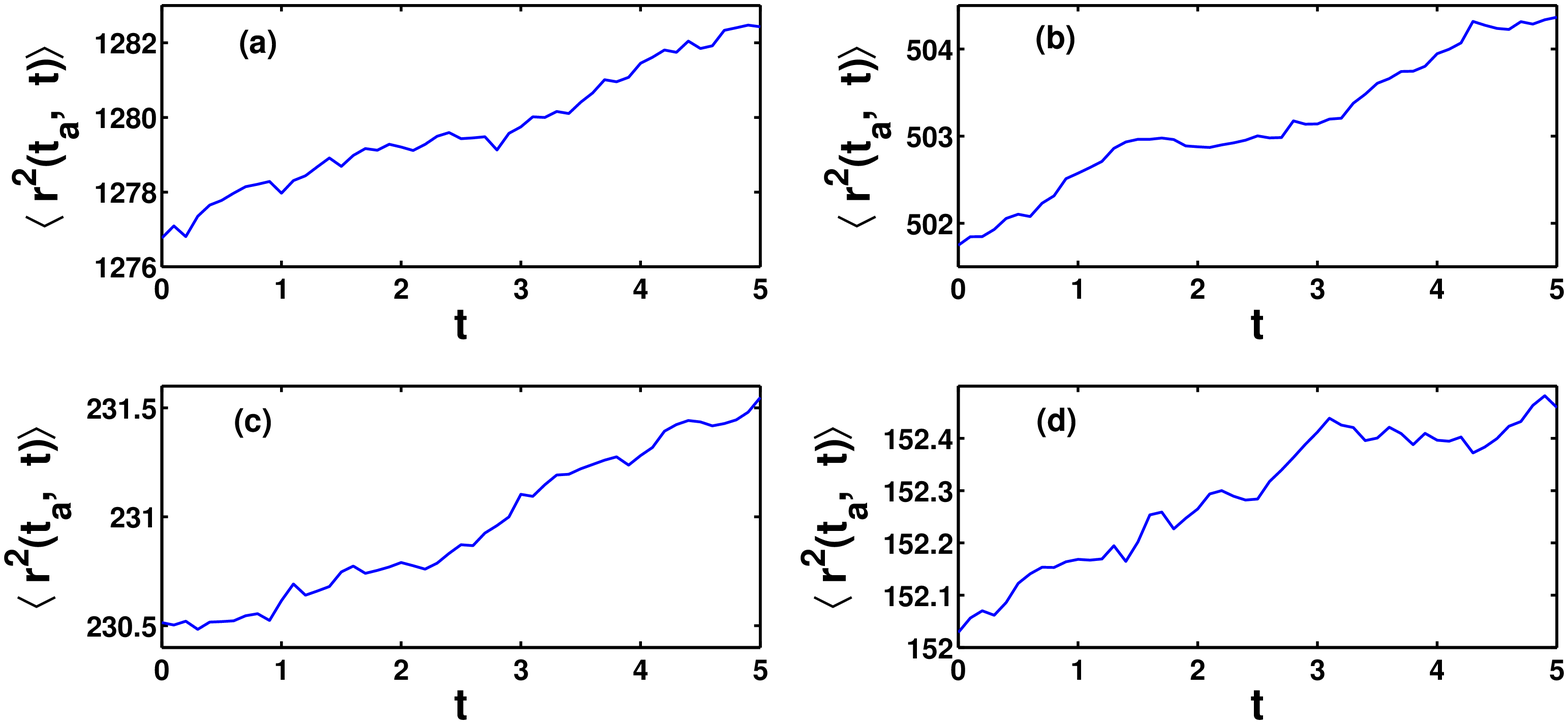}\\
  \caption{The relation between $\langle r^2(t_{a},  t)\rangle $ and the observation time $t$ for $t_{a} \gg t$.  It can be noted  that $\langle r^2(t_{a},  t)\rangle $ increases linearly with $t$ for fixed $t_{a}$ and the fluctuations are large because of population splitting.  The parameter $t_{a}=1000$,   $t=5$,   $\lambda=10^{-1}$ for (a), $\lambda=10^{-2}$ for (b), $\lambda=10^{-3}$ for (c), $\lambda=10^{-4}$ for (d),  $\alpha=0.6$  and the number of the trajectories is $10000$.   }\label{agfig22}
\end{figure}

\begin{figure}[htb]
  \centering
  \includegraphics[width=9cm,  height=6cm]{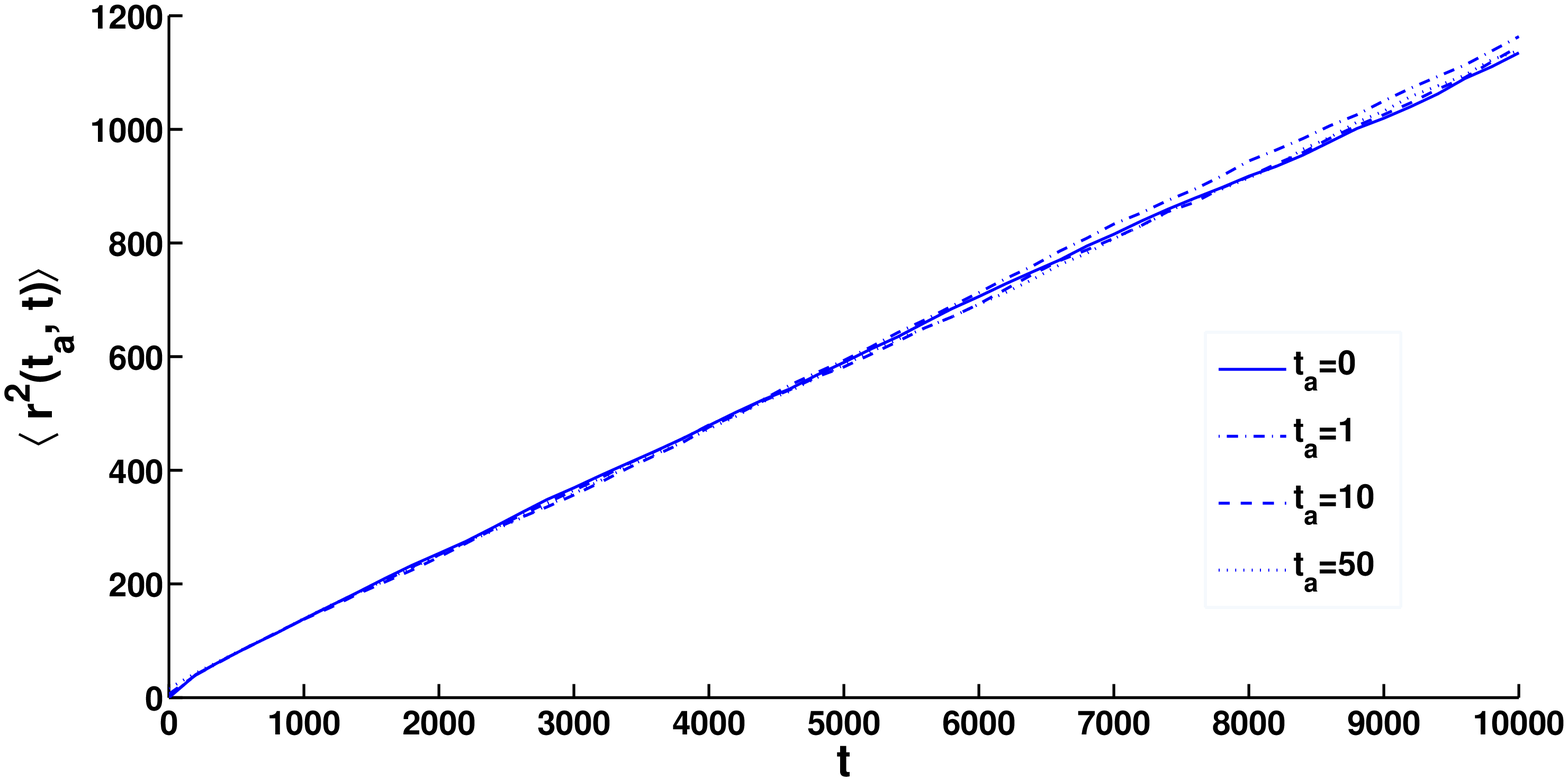}\\
  \caption{The relation between $\langle r^2(t_{a},  t)\rangle $ and the observation time $t$ for various $t_a$ with $t \gg t_a$.  }\label{agfig23}
\end{figure}

\subsection{Propagator function  $P(x, t_{a}, t)$}\label{subsec3-2}
\begin{figure}[htb]
  \centering
  \includegraphics[width=9cm,  height=6cm]{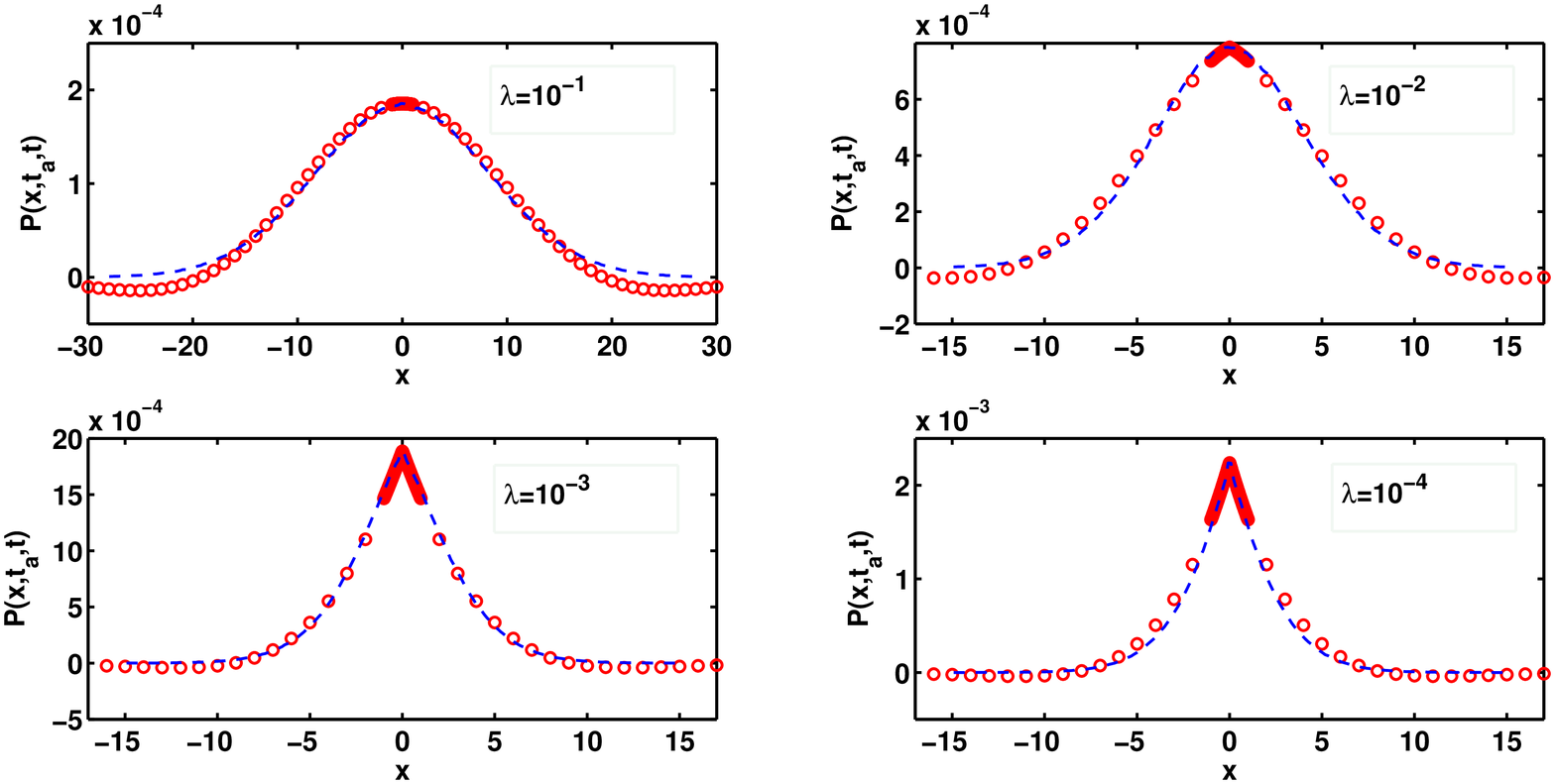}
  \caption{Propagator functions with $t_{a}=3$, $\alpha=0.6$, and $t=500$. The red lines with empty symbol are obtained by calculating Eq. (\ref{ageq42}), and the blue dash lines are got from the generated $10^4$ trajectories of the particles.   }\label{agfig51}
\end{figure}

 %We consider the simulation of $P(x,  t_{a},  t)$.  we don't obtain the accurate form of $P(x, t_{a}, t)$,   In order to have a better understand of the propagator function of the tempered ACTRW,   we use numerical inversion of Laplace transform.
In this subsection, we discuss the propagator function $P(x, t_{a}, t)$ of the tempered ACTRW. Omitting the motionless part of Eq. (\ref{ageq303}), taking $f(x)$ as Gaussian, and performing inverse Fourier transform w.r.t. $k$, there exists
\begin{equation}\label{ageq41}
  P(x,  s,  u) \sim \frac{\omega(s,  u)}{2u}F_{1}(u,x)
\end{equation}
with
\begin{equation*}
\begin{split}
  F_{1}(u,x) = & \sqrt{\frac{(u+\lambda)^\alpha-\lambda^\alpha}{0. 5(1+\lambda^\alpha-(u+\lambda)^\alpha)}} \\
    & \cdot \exp \left(-|x|\sqrt{\frac{(u+\lambda)^\alpha-\lambda^\alpha}{0. 5(1+\lambda^\alpha-(u+\lambda)^\alpha)}}\right).
\end{split}
\end{equation*}
For $u \ll s$, Eq. (\ref{ageq41}) can be rewritten as
\begin{equation}\label{ageq42}
  P(x,  t_{a},  u)   \sim \left(\frac{\lambda^\alpha-(u+\lambda)^\alpha}{2u}[g(t_a)*1]+1\right)F_1(u,x).
\end{equation}
From FIG. (\ref{agfig51}), it can be noted that for small $\lambda$ ($\lambda=10^{-3}$ or $\lambda=10^{-4}$) the propagator functions display the characteristics of $\alpha$ stable distribution; while for large $\lambda$ the $P(x,t_{a},t)$ shows the classical normal behavior.

%\begin{widetext}
%\begin{equation}\label{ageq52}
%P(x,  t_{a},  u)=\left\{
 %              \begin{array}{ll}
 %               \frac{\lambda^{\alpha}-(u+\lambda)^{\alpha}}{2u}\int_{0}^{t_{a}}\tau^{\alpha-1}exp(-\lambda \tau)E_{\alpha,  \alpha}(\lambda^{\alpha}\tau^{\alpha})d\tau M_{3}(u,  x),   & \hbox{ for $t_{a}\ll t$;} \\
  %              \frac{(u+\lambda)^{\alpha}-\lambda^\alpha}{2u^2}t_{a}^{\alpha-1}exp(-\lambda t_{a})E_{\alpha,  \alpha}(\lambda^{\alpha}t_{a}^{\alpha})M_{3}(u,  x) ,   & \hbox{~for $t_{a}\gg t$. }
   %            \end{array}
    %         \right.
%\end{equation}
%\end{widetext}
%\begin{equation}\label{ageq42}
%  P(x,  t_{a},  u) \sim \frac{1}{2u}F_1(u,x).
%\end{equation}

\begin{figure}[htb]
  \centering
  \includegraphics[width=9cm,  height=6cm]{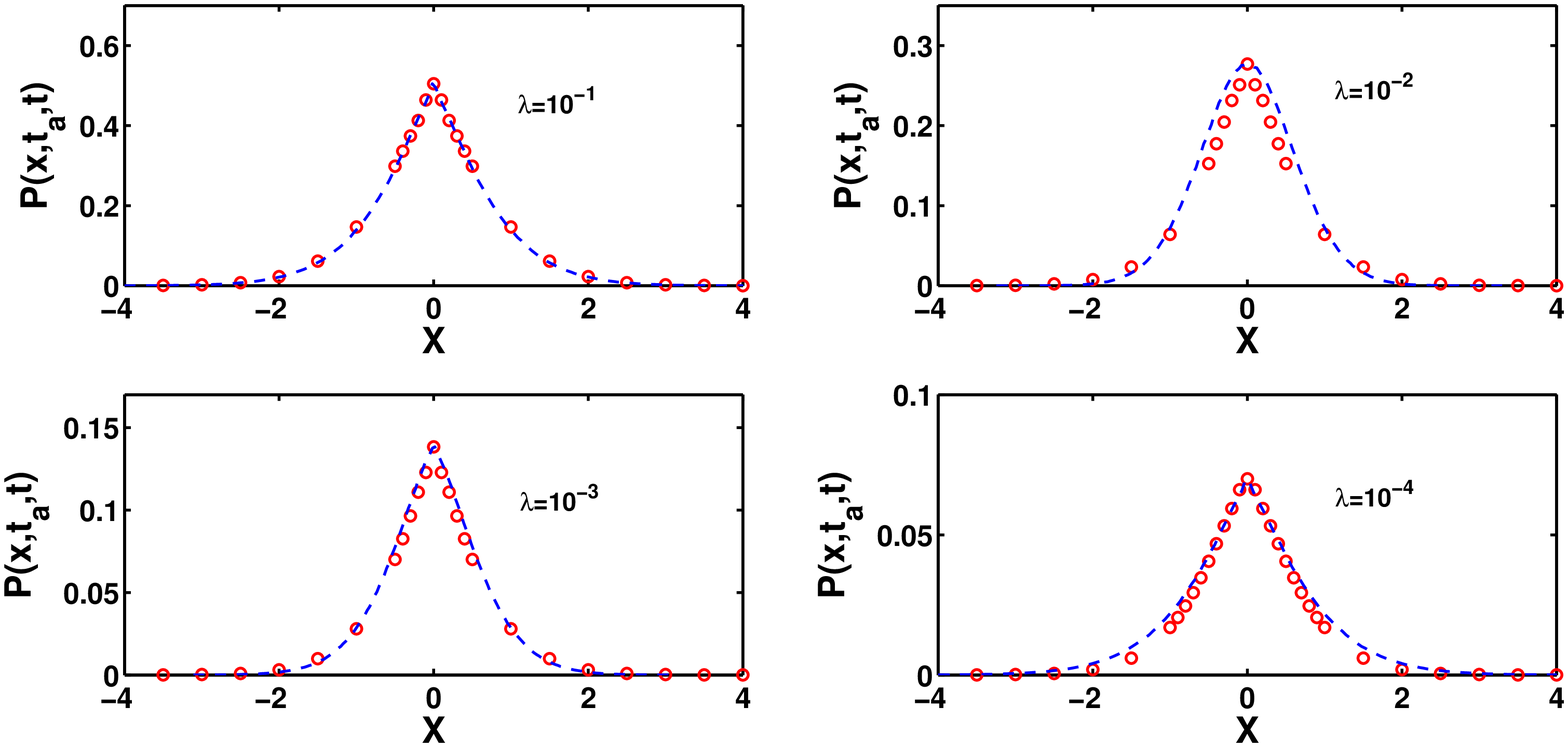}\\
  \caption{Strong aging case (contrary to FIG. (\ref{agfig51})) with $t_{a}=500$ and $t=3$. The other parameters are same as FIG. (\ref{agfig51}). }\label{agfig52}
\end{figure}
%\begin{figure}[htb]
%  \centering
%  % Requires \usepackage{graphicx}
%  \includegraphics[width=9cm,  height=6cm]{agfig53.eps}\\
%  \caption{Relation between $P(x,  t_{a},  t)$ and $t_a$ for strong aging from Eq.(\ref{ageq43}).  %from the figure we can know that for large $\lambda$ the process tend to equilibrium easily.
%  }\label{agfig53}
%\end{figure}
For $u \gg s$, Eq. (\ref{ageq303}) yields,
\begin{equation}\label{ageq43}
  P(x,  t_{a},  u)   \sim \frac{(u+\lambda)^\alpha-\lambda^\alpha}{2u^2}g(t_a)F_1(u,x).
\end{equation}
Contrary to  FIG. (\ref{agfig51}), FIG. (\ref{agfig52}) displays the behaviors of the $\alpha$ stable distribution for all kinds of $\lambda$.

From the numerical results and the theory we can see that the `$\alpha$ stable distribution' characteristics can be found for small $\lambda$. Both the distributions for $\lambda=10^{-3}$ and $\lambda=10^{-4}$ have the sharp peak and the tail of Eq. (\ref{ageq42}) decays slowly. While for large $\lambda$, the top of the distribution for  Eq. (\ref{ageq42}) is smooth, being different from the case of small $\lambda$. Therefore, depending on the choice of $\lambda$, one can control the behaviors of the propagator.

\section{strong relation between the fluctuation and response}\label{sec4}
In this section, we discuss the aging from a new point of view. Based on the CTRW model, consider such a process: the particles begin to move at time $t=0$ and undergo unbiased diffusion in the time interval $(0,t_a)$; then an external field is switched on the system starting from $t_a$. If the averaged response of the particles depends on $t_a$, the process is said to exhibit aging. Generally speaking, giving some disturbance to a system, some characteristics (parameters of thermodynamics) of the system will change, being called response \cite{Bertin:1}. Under the small disturbance of external field, if the change of the parameter of thermodynamics is proportional to the force of external field, then it is called linear response. It seems important to use drift diffusion to consider aging. Using the method given in \cite{Barkai:7,Allegrini:3,Froemberg:1,Shemer:1}, we discuss the tempered aging Einstein relation.

%In this parts, we consider the aging in a new view, Based on the model of CTRW. We consider such a process: the particles begin  to move at time $t=0$,
%the particle move without deviation in the time between $0$ and $t_{a}$, contrary to the symmetry case, we give a force $F$ at time $t_{a}$, which destroy the symmetrical characteristic.
%when the ensemble average is rely on $t_{a}$, such a process show aging. Generally speaking, Giving some disturbance to the system, then some character (parameters of thermodynamics)of the system will change, we call this response \cite{Bertin:1}. If external field (disturbance) is small, then the parameters of thermodynamics is  in proportion to external field, which is called linear response. Therefore it's of importance to use drift diffusion to consider aging. Using the methods\cite{Barkai:7,Allegrini:3,Froemberg:1,Shemer:1} give the templed Einstein relation .

Let us consider a simple example of random walk on a one-dimensional lattice; the length of the lattice is $c$, and the particles can only move to its neighboring sites. Waiting times of different steps of the random walk are considered independent and have the same distribution $\psi(t)$. Jumps to the right (left) are performed with the probability $1/2+h/2$ ($1/2-h/2$). The total time is $t=t_{a}+t_{b}$; $[0, t_{a}]$ is called aging interval with $h=0$; and $(t_{a},t_{a}+t_{b})$ is called response interval with $0<h<1$. Let $x=x_{a}+x_{b}$, where $x_{a}=\sum_{i=1}^{n_{a}}x_{i}^{a}$ is the displacement performed in the aging time  interval and $x_{b}=\sum_{i=1}^{n_{b}}x_{i}^{b}$ is the displacement performed in the response time interval, and  $x_{i}^{a}, x_{i}^{b}$ are the step lengths and $n_{a}, n_{b}$ are the number of events happened in the two time intervals, repectively.

We consider the correlation function $\langle  (x_{a})^2x_{b}\rangle$ which shows the impact between $(x_{a})^2$ in the aging interval and $x_{b}$ in the response interval. And define a parameter $F_{R}$ to show the relation between fluctuation and response \cite{Barkai:7},
\begin{equation}\label{ageq601}
F_{R}=\frac{\langle (x_{a})^2x_{b}\rangle}{\langle (x_{a})^2\rangle \langle x_{b} \rangle}-1.
\end{equation}
If $F_{R}=0$, it shows that $x_{a}^2$ and $x_{b}$ are independent with each other.
Using the relation $\langle x_{a}^2 \rangle=c^2\langle n_{a}\rangle$ and $\langle x_{b}\rangle =hc\langle n_{b}\rangle$, then $F_{R}$ can be shown in another way,
\begin{equation}\label{ageq602}
F_{R}=\frac{\langle n_{a}n_{b}\rangle}{\langle n_{a}\rangle \langle n_{b} \rangle}-1.
\end{equation}
We further introduce $X_{t_{a},t_{b}}(n_{a},n_{b})$, the probability to occur $n_{a}$ events in the aging interval and $n_{b}$ events in the response interval. Following the result given in \cite{Luck:1},
\begin{equation}\label{ageq603}
\begin{split}
   X_{t_{a},t_{b}}(n_{a},n_{b}) & =\langle I(t_{n_{a}}<t_{a}<t_{n_{a}+1}) \\
    & I(t_{n_{a}+n_{b}}<t_{a}+t_{b}<t_{n_{a}+n_{b}+1})\rangle,
\end{split}
\end{equation}
where $I(t_{n_{a}}<t_{a}<t_{n_{a}+1})=1$ if the event inside the parenthesis occurs, and $0$ if not.
Using double Laplace transform, if $n_b=0$,
\begin{equation}\label{ageq604}
    X_{s,u}(n_{a},n_{b}) =\frac{\varphi^{n_{a}}(s)}{u}\left[\frac{1-\varphi(s)}{s}-\frac{\varphi(u)-\varphi(s)}{s-u}\right];
\end{equation}
and if $n_b\geq 1$,
\begin{equation*}
 X_{s,u}(n_{a},n_{b}) =\frac{\varphi^{n_{a}(s)}\psi^{n_{b}-1}(u)}{u(s-u)}[1-\varphi(u)][\varphi(u)-\varphi(s)].
\end{equation*}
Summing $n_a$, $n_b$ from $0$ to $\infty$  leads to
\begin{equation*}
  \langle n_a n_b  \rangle_{s,u}=\frac{[\varphi(u)-\varphi(s)]\varphi(s)}{u(s-u)(1-\varphi(s))^2(1-\varphi(u))}.
\end{equation*}
Using the Laplace transform of $\varphi(t)$ (\ref{ageq12}), we have, when $t_{a}\gg t_{b}$,
\begin{equation}\label{ageq6014}
\begin{split}
  \langle n_{a}n_{b} \rangle_{t_{a},t_{b}} \sim  & \, t_{b}[t_{a}^{\alpha-1}E_{\alpha,\alpha}(\lambda^\alpha t_{a}^{\alpha})\exp(-\lambda t_{a})] \\
    & *[t_{a}^{\alpha-1}E_{\alpha,\alpha}(\lambda^\alpha t_{a}^{\alpha})\exp(-\lambda t_{a})];
\end{split}
\end{equation}
if $t_{a} \ll t_{b}$, there exists
\begin{equation}\label{ageq6012}
\begin{split}
\langle n_{a}n_{b} \rangle_{t_{a},t_{b}}  \sim & (1*[t_{a}^{\alpha-1}E_{\alpha,\alpha}(\lambda^\alpha t_{a}^{\alpha})\exp(-\lambda t_{a})]) \\
    & \cdot(1*[t_{b}^{\alpha-1}E_{\alpha,\alpha}(\lambda^\alpha t_{b}^{\alpha})\exp(-\lambda t_{b})])\\
    &-1*[t_{a}^{\alpha-1}E_{\alpha,\alpha}(\lambda^\alpha t_{a}^{\alpha})\exp(-\lambda t_{a})]\\
   &*[t_{a}^{\alpha-1}E_{\alpha,\alpha}(\lambda^\alpha t_{a}^{\alpha})\exp(-\lambda t_{a})] .
\end{split}
\end{equation}
Following the results given by \cite{Luck:1},
\begin{equation*}
  \langle n_a\rangle_{s,u}=\frac{\varphi(s)}{us(1-\varphi(s))},
\end{equation*}
there exists
\begin{equation}\label{ageq605}
\langle n_{a} \rangle_{t_{a},t_{b}}=(1*[t_{a}^{\alpha-1}E_{\alpha,\alpha}(\lambda^\alpha t_{a}^{\alpha})\exp(-\lambda t_{a})]).
\end{equation}
And $\langle n_b\rangle_{s,u} $ is the same as Eq. (\ref{ageq110}).
%\begin{equation}\label{ageq604}
%  X_{s,u}(n_{a},n_{b})=\left\{
%                         \begin{array}{ll}
%                           \frac{\psi^{n_{a}}(s)}{u}[\frac{1-\psi}{s}-\frac{\psi(u)-\psi(s)}{s-u}], & \hbox{$n_{b}=0$;} \\
%                           \frac{\psi^{n_{a}(s)\psi^{n_{b}-1}(u)}}{u(s-u)}[1-\psi(u)][\psi(u)-\psi(s)], & \hbox{else.}
%                         \end{array}
%                       \right.
%\end{equation}
%for $t_{a}>>t_{b}$ ,$ \langle n_{b} \rangle_{t_{a},t_{b}}$ can be shown by
%\begin{equation}\label{ageq606}
% \langle n_{b} \rangle_{t_{a},t_{b}}\sim t_b[t_{a}^{\alpha-1}E_{\alpha,\alpha}(\lambda^\alpha t_{a}^{\alpha})exp(-\lambda t_{a})]
%\end{equation}
%while for $t_{a}<<t_{b}$  yields,
%\begin{equation}\label{ageq6013}
% \langle n_{b} \rangle_{t_{a},t_{b}}  \sim (t_{b}^{\alpha-1}E_{\alpha,\alpha}(\lambda^\alpha t_{b}^{\alpha})exp(-\lambda t_{b}))
%\end{equation}
For $t_a \gg t_b$, there exists
\begin{widetext}
\begin{equation}\label{ageq607}
F_{R}(t_{a},t_{b}) \sim \frac{[t_{a}^{\alpha-1}E_{\alpha,\alpha}(\lambda^\alpha t_{a}^{\alpha})\exp(-\lambda t_{a})]*[t_{a}^{\alpha-1}E_{\alpha,\alpha}(\lambda^\alpha t_{a}^{\alpha})\exp(-\lambda t_{a})]}{(1*[t_{a}^{\alpha-1}E_{\alpha,\alpha}(\lambda^\alpha t_{a}^{\alpha})\exp(-\lambda t_{a})])t_{a}^{\alpha-1}E_{\alpha,\alpha}(\lambda^\alpha t_{a}^{\alpha})\exp(-\lambda t_{a})}.
\end{equation}
\end{widetext}
From Eq. (\ref{ageq607}), we know that $F_{R}(t_{a},t_{b})$ does depend on $t_{b}$.

In the following, we further consider the Einstein relation \cite{Froemberg:1,Shemer:1} for the tempered aging process.
Denoting $\langle x(t_{a},t_{b})\rangle_{F}$ as the first order moment of the displacement under the influence of a force $F$, from Eq. (\ref{ageq114}) and $\langle x(t_{a},t_{b})\rangle_{F}=hc \langle n_b \rangle$, we get that for $t_{a}\gg t_b$,
%the first order moment of the displacement under the influence of a force, denoted by $\langle x(t_{a},t_{b})\rangle_{F}$, is given as
% we can easy to get the first moments
%of displacement, we denote it $\langle x(t_{a},t_{b})\rangle_{F}$, and the system  in under the influence of a force.
\begin{equation}\label{ageq609}
\langle x(t_{a},t_{b})\rangle_{F} \sim hc t_{b}t_{a}^{\alpha-1}E_{\alpha,\alpha}(\lambda^\alpha t^\alpha_{a})\exp(-\lambda t_{a}).
\end{equation}
Denoting $\langle r^2(t_{a},t_{b})\rangle_0$ as the mean squared displacement of the random walk without external force, from Eq. (\ref{ageq306}) we obtain that for $t_{a} \gg t_b$,
\begin{equation}\label{ageq6010}
  \langle r^2(t_{a},t_{b})\rangle_{0} \sim M_{2}t_{b}t_{a}^{\alpha-1}E_{\alpha,\alpha}(\lambda^\alpha t^\alpha_{a})\exp(-\lambda t_{a})
\end{equation}
with $M_{2}=1/2c^2+1/2(-c)^2=c^2$. Under the assumption $h=\frac{cF}{2K_{b}T} \ll 1$, we obtain the following relation,
\begin{equation}\label{ageq6011}
 \langle x(t_{a},t_{b})\rangle_F \sim \frac{F}{2K_{b}T} \langle r^2(t_{a},t_{b})\rangle_{0}.
\end{equation}
%the results is same to $\phi(t) \sim t^{-\alpha-1}$. From Eq. (\ref{ageq6011}), when we know each of both, we can obtain the other.

\section{Fokker-Planck equation for the tempered ACTRW}\label{sec5}
We now derive the Fokker-Planck equation of the tempered ACTRW, which can be used to solve the tempered aging diffusion problems with different types of boundary and initial condition.
%We now consider the aging diffusion equation which describes a large number of ACTRW process. The equation is characterized by fractional order, which is called as fractional order diffusion equation\cite{Barkai:6,Sokolov:1,Schulz:2}. It can be used to solve aging diffusion problems with different types of boundary and initial conditions. we call this equation is the Fokker-Planck equation for ACTRW\cite{Barkai:2}.
Omitting the motionless part of Eq. (\ref{ageq303}) and taking $f(x)$ as Gaussian (i.e., $f(k) \sim 1- \frac{1}{2}k^2$), we have
%\begin{widetext}
\begin{equation}\label{ageq51}
\begin{array}{l}
P(k,  s,  u) \\
 \displaystyle\sim \frac{\omega(s,  u)}{u}\frac{(u+\lambda)^\alpha-\lambda^\alpha}{(u+\lambda)^\alpha-\lambda^\alpha
+\frac{1}{2}k^2(1+\lambda^\alpha-(u+\lambda)^\alpha)}.
\end{array}
\end{equation}
%\end{widetext}
Define the Riemann-Liouville fractional derivative from $0$ to $t$ as,
\begin{equation}\label{ageq52}
_{0}D_{t}^{p}y(t)=\frac{1}{\Gamma(n-p)}\frac{d^{n}}{d t^n }\int_{0}^{t}(t-\tau)^{n-p-1}y(\tau) d \tau,
\end{equation}
where $n=[p]+1$ and $[\, ]$ denotes the integer part of $p$.
For $0<p<1$, taking the Laplace transform of the Riemann-Liouville fractional derivative results in
\begin{equation}\label{ageq53}
\mathcal{L} [ {_{0}D_{t}^{p}y(t)} ]=s^p \mathcal{L} [y(t)]=s^p y(s).
\end{equation}
From the property of the inverse Fourier transform, we have
\begin{equation}\label{ageq54}
\mathcal{F}^{-1}[k^2y(k)]=-\frac{\partial^2}{\partial x^2}y(x).
\end{equation}
From Eq. (\ref{ageq53}) and Eq. (\ref{ageq54}), performing once inverse Fourier transform and double inverse Laplace transforms leads to
%Now,  we calculate the Fokker-Planck equation for the tempered ACTRW, using Eq.(\ref{ageq53}) and (\ref{ageq54}), Performing once inverse fourier transform and double inverse laplace of Eq.(\ref{ageq51}),  we can derive,  \\
\begin{widetext}
\begin{equation}\label{ageq55}
\begin{split}
&-e^{-\lambda t}~_{0}D_{t}^{\alpha}\left(e^{\lambda t}P(x,  t_a,  t)\right)+\lambda^\alpha P(x,t_{a},t)\\
=&-\frac{1}{2}(1+\lambda^\alpha)\frac{\partial^2}{\partial x^2}  P(x,  t_a,  t)
+e^{-\lambda t}{_{0}D_{t}^{\alpha}} \left( e^{\lambda t}\frac{\partial^2}{\partial x^2}  P(x,  t_a,  t) \right)\\
 &+\lambda^{\alpha}[\omega(t_{a},t)\ast 1]\delta(x)-[1 \ast e^{-\lambda t} {_{0}D_{t}^{\alpha}}\left(e^{\lambda t}\omega( t_a,  t)\right)]\delta(x),
\end{split}
\end{equation}
\end{widetext}
where the notation $\ast$ represents the convolution of the functions w.r.t. $t$. Eq. (\ref{ageq55}) is the Fokker-Planck  equation of the Green function $P(x, t_a, t)$ in the case that the waiting time distribution is the tempered power-law (\ref{ageq11}).

For the non tempered case, namely, $\lambda=0$ and $\varphi(t)\sim t^{-1-\alpha}$,  %being a pure power-law distribution. when $0<\alpha<1$,  there are also many properties in aging systems and we can see obvious aging effects. In our follows,  we give the aging diffusion equation in this case,
taking $\lambda=0$ in Eq. (\ref{ageq51}) results in
\begin{equation}\label{ageq56}
P(k,  s,  u)\sim \omega(s,  u)\frac{u^{\alpha-1}}{u^{\alpha}+\frac{k^{2}}{2}(1-u^\alpha )}.
\end{equation}
Performing inverse Fourier transform and double inverse Laplace transform on Eq. (\ref{ageq53})  yields the corresponding aging diffusion equation
%\begin{widetext}
\begin{equation}\label{ageq57}
\begin{split}
_{0}D_{t}^{\alpha}P(x,  t_a,  t)=&\frac{1}{2}\frac{\partial^2}{\partial x^2}P(x,t_{a},t)-{_{0}D}_{t}^{\alpha}\frac{\partial^2}{\partial x^2}P(x,t_{a},t)
\\
& -\omega(t_a,  t)\ast\frac{t^{-\alpha}}{\Gamma(1-\alpha)}\delta(x).
\end{split}
\end{equation}
%\end{widetext}
As expected, taking $\lambda=0$ in Eq. (\ref{ageq55}) also arrives at Eq. (\ref{ageq57}). There are also other forms of Eq. (\ref{ageq55}), e.g., adding the motionless part of (\ref{ageq303}) to the equation; for the longer time scale $t \gg 1$, then $u^\alpha k^2$ can be reasonably omitted \cite{Barkai:5}, i.e., the term $_{0}D_{t}^{\alpha}\frac{\partial^2}{\partial x^2}P(x,t_{a},t)$ can be omitted in Eq. (\ref{ageq57}).

%Comparing Eq. (\ref{ageq55}) to Eq.(\ref{ageq57}), we can see if make $\lambda=0$ in Eq.(\ref{ageq55}),  we can also get the Eq .(\ref{ageq57}),  therefore,  eq.(\ref{ageq57}) is a limit form of eq. ((\ref{ageq55})).  Of course,  there are many other forms of Eq.(\ref{ageq55}).  For example,  we can consider both the moving and the not moving part of $\ref{ageq303}$, From this case, we can derive the different diffusion equations of ACTRW.
%Then whatever the case are,   they can describe a large number of ACTRW process.
%For longer time scale $t>>1$, then $ u^\alpha<<1 $, then $u^\alpha k^2$ can be omit, i.e. the term $_{0}D_{t}^{\alpha}\frac{\partial^2}{\partial x^2}P(x,t_{a},t)$ can be omit in Eq.(\ref{ageq57})

\section{Conclusions}
Because of the boundedness of physical space and the finiteness of the lifetime of particles, sometimes it is a more physical choice to use tempered power-law jump length or waiting time distribution instead of the pure power-law distribution. This paper discusses the renewal process and ACTRW with the tempered power-law waiting time distribution $\varphi(t)\sim e^{-\lambda t}t^{-1-\alpha}$. Since the tempered power-law distribution lies between the pure power-law and exponential distributions, as expected, the transition dynamics is found with the time evolution and the turning point depends on $\lambda$. By using the aging renewal theory, the $p$-th moments of the renewal times $n_a(t_a, t)$ are analytically obtained and numerically confirmed by simulating the particles' trajectories. In particular, the first order moment of $n_a(t_a, t)$ is more detailedly discussed.   Similarly, the mean squared displacement of the tempered ACTRW is analytically got and numerically verified. Based on the $p$-th moment of $n_a(t_a, t)$ and mean squared displacement of the tempered ACTRW, the aging effects are deeply analyzed. Finally, the tempered aging Einstein relation is attained and the corresponding Fokker-Planck equation is derived from the tempered ACTRW.

%The aging effect of the
%
%In this paper,  we propose a tempered distribution of waiting time,  $\varphi(t)\sim e^{-\lambda t}t^{-1-\alpha}$. And using this distribution to discuss the aging effect. We can see when $t\rightarrow \infty$,  the PDF of the waiting time is similar to the Poisson distribution. In this case,  we can not observe the aging phenomenon.
%We use the aging renewal theory to obtain the pth moment of the renewal process $n_a(t_a, t)$, and discuss the properties of $\langle n_a(t_a, t)\rangle$. especially, we analyze the characters of the survival probability in slightly and strongly aged systems, respectively. We then turn to describe the tempered ACTRW. In this case of waiting time distribution, we obtain the mean square displacement in two limit conditions, and we can demonstrate these effectiveness. In Sec.(\ref{sec4}) we discuss the propagator by the numerical inversion of Laplace transform. For the aging diffusion equation, we obtain the Fokker-Planck equation for the tempered ACTRW , which is a fractional diffusion equation and can describe a large class of aging phenomenons.

\section*{Acknowledgments}
The authors thank Eli Barkai for the discussions. This work was supported by the Fundamental Research Funds for the Central Universities under Grant No. lzujbky-2015-77, and the National Natural Science Foundation of China under Grant No. 11271173.

%\section*{Appendix A} NO NEED FOR APPENDIX.  WHEN YOU DISCUSS THE
%LANGEVIN SIMULATIONS D A SENTENCE.  To obtain trajectories of
%fractional Langevin s.  we extended the Hosking method with ideas
%from \cite{Deng:07}.  AT IS ENOUGH.
% The Hosking method \cite{Hosking:84} is used in our
%simulation of fractional Brownian motion,   which is an exact method.
%This method first generates the sample of fractional Gaussian noise,
%then by taking the cumulative sum for the noise samples the
%fractional Brownian motions sample is obtained.  Here it should be
%emphasized that because of the `historical dependence' (long-range
%correlations) of fractional Gaussian noise when we compute the
%(n+1)th noise,   all the information of the first $n$ noise samples
%must be used.  So the computational cost is greater compare with the
%simulation of CTRW model being completely independent for each step.
%The numerical algorithm \cite{Numerical:08}
% for simulating the trajectories of
%generalized Langevin equation is developed by combining the idea
%(predictor-corrector) in our previous work \cite{Deng:07} with
%Hosking method.
\appendix
%\appendixpage
%\addcontentsline{toc}{chapter}{¸½Â¼}
%\markboth{¸½Â¼}{}
\begin{appendices}
\section{Generation of random variables (FIG. 1)} \label{AppendixC}

When generating the random variables with the PDF Eq. (\ref{ageq11}) to plot FIG. \ref{agfig11}, the Monte Carlo statistical methods \cite{Christain:1} is used. We first rewrite $\varphi(t)$ as $\varphi(t)=H(x)f_1(x)$, where $f_1(x)=\alpha t_{0}^\alpha t^{-\alpha-1}$ with $t_0$ being a small number. Denote the maximum of $H(x)$ as $M$. Then the algorithm can be described as:
%\textcolor[rgb]{1.00,0.00,0.00}{In numerous setting, the distribution for the density $\varphi(t)$ is difficult to  simulated directly because of the complexity of function $\varphi(t)$, which may require substantial compute time or the the original function doesn't exit.  We use Monte Carlo statistical methods \cite{Christain:1} to simulated  FIG. 1,
%if f can spread out into two function, i.e. $f(X)=H(x)f_{1}(x)$, here $f_{1}(x)$ is a arbitrary density function and $H(x)$ has a maximum $M$ for every $x$. the method as follows,
\begin{enumerate}
  \item Generate a r.v. $x_{f_{1}}$ with PDF $f_{1}$ and a r.v. $\xi$ being uniformly distributed in the interval $[0,1]$.
  \item Accept $x_{f_{1}}$, if $M\xi\leq H(x_{f_{1}})$; otherwise, reject.
  \item Return to Step $1$.
\end{enumerate}
% in our simulation $f_{1}(x)=\alpha t_{0}^\alpha t^{-\alpha-1}$, here $t_{0}$ is a small number.
%}

%\begin{equation}\label{ageq931}
%\left\{
%  \begin{array}{ll}
%    1, & \hbox{Generate a random variables $x_{f_{1}}$,
%which has a $f_{1}$ distribution  and $\xi \sim U[0,1]$;} \\
%    2, & \hbox{Accept $X_{f_{1}}$, if $M\xi\leq H(x_{f_{1}})$,then continue;} \\
%    3, & \hbox{Otherwise, reject, continue.}
%  \end{array}
%\right.
%\end{equation}

\section{Mittag-Leffler function} \label{AppendiA}
The two-parameter function of the Mittag-Leffler type plays a very important role in the
fractional calculus, being introduced by G.M. Mittag-Leffler and studied
by A. Wiman \cite{podlubny:1}. %A lot of relationships for this function were obtained by Humbert and Agarwal using the Laplace
%transform technique, and this function also called Agarwal function.
The two-parameter Mittag-Leffler function is defined by the series expansion
\begin{equation}\label{ageq911}
 E_{\alpha,\beta}=\sum_{k=0}^{\infty} \frac{z^{k}}{\Gamma(\alpha k+\beta)}
\end{equation}
with $\alpha>0$ and $\beta>0$. Its one-parameter form ($\beta=1$) is given as
\begin{equation}\label{ageq912}
 E_{\alpha}=\sum_{k=0}^{\infty} \frac{z^{k}}{\Gamma(\alpha k+1)}.
\end{equation}
The asymptotic expansions of Mittag-Leffler are important for obtaining the various useful estimates of the long time or short time fractional dynamics. For small $z$, there exists
\begin{equation}\label{age913}
 E_{\alpha,\beta}(z) \sim \frac{1}{\Gamma(\beta)}+\frac{z}{\Gamma(\alpha+\beta)},
\end{equation}
in the special case, $ E_{\alpha}(z) \sim 1+\frac{z}{\Gamma(\alpha+1)}$.
Another important and useful formula is the asymptotic expansion of large scale for the Mittag-Leffler function. For $0<\alpha<2$, $\beta$ is an arbitrary complex number and $\mu$ is an arbitrary real number such that $\pi\alpha/2< \mu< \min\{ \pi, \pi\alpha\}$, then for an arbitrary integer $p\geq 1$, the following expansion holds,
\begin{equation}\label{ageq914}
 \begin{split}
   E_{\alpha,\beta}(z) \sim & \frac{1}{\alpha}z^{(1-\beta)/\alpha}\exp(z^{1/\alpha}) \\
     & -\sum_{k=1}^{p}\frac{z^{-k}}{\Gamma(\beta-\alpha k)}
 \end{split}
\end{equation}
with $z\rightarrow \infty$ and $\arg(z) \leq \mu$. For $z\rightarrow +\infty$, from Eq. (\ref{ageq914}), we have
\begin{equation}\label{ageq915}
E_{\alpha,\beta}(z) \sim \frac{1}{\alpha}z^{(1-\beta)/\alpha}\exp(z^{1/\alpha}).
\end{equation}
And if $z\rightarrow -\infty$, %the first part convergent to $0$ fastly,
\begin{equation}\label{age916}
 E_{\alpha,\beta}(z)\sim -\frac{1}{z \Gamma(\beta-\alpha)}-\frac{1}{z^2 \Gamma(\beta-2\alpha)};
\end{equation}
when $\alpha=\beta$, we have $|\Gamma(0)|\rightarrow \infty $. Therefore  $ E_{\alpha,\alpha}(z)\sim -\frac{1}{z^2 \Gamma(-\alpha)}$.

\begin{figure}[htb]
  \centering
  \includegraphics[width=9cm, height=6cm]{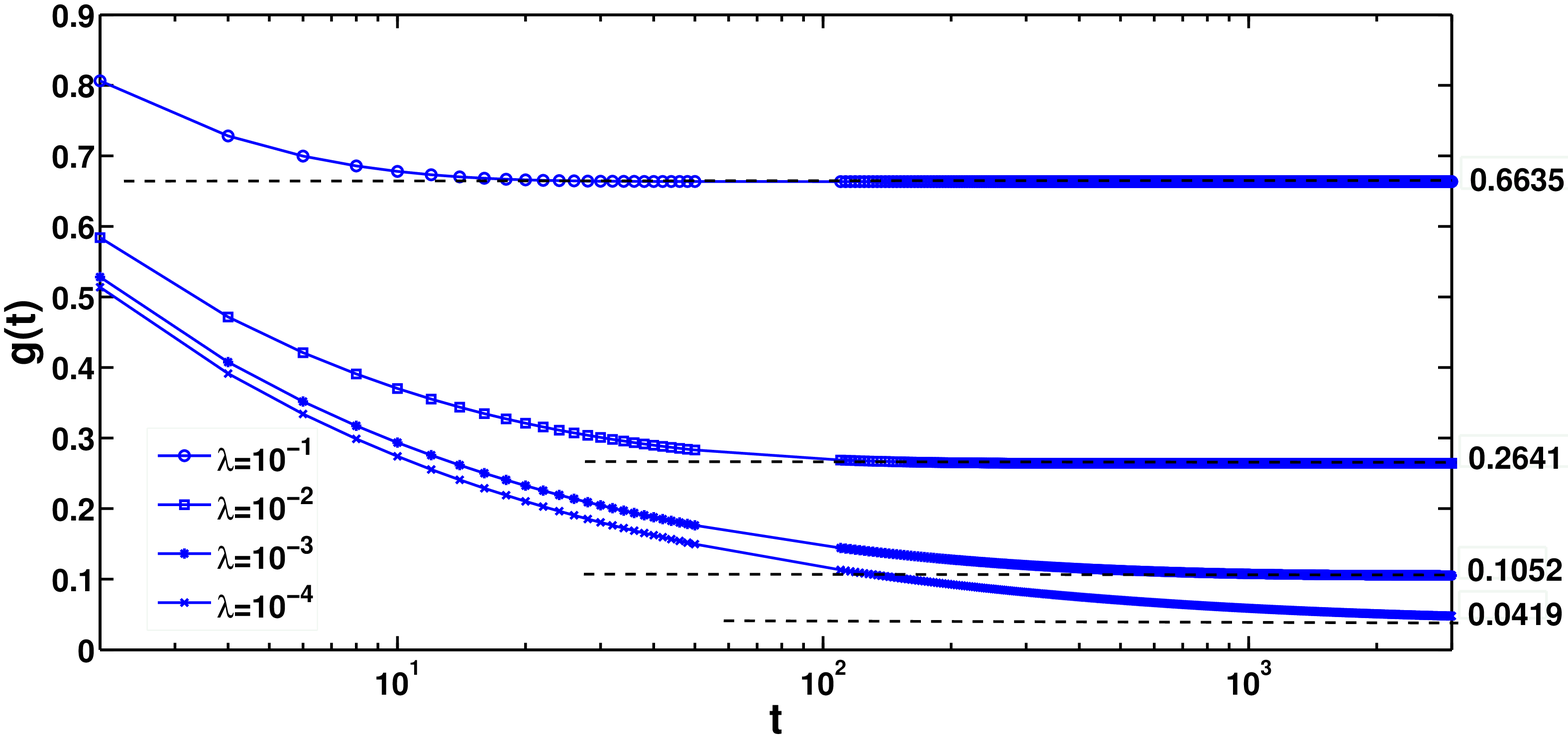}\\
  \caption{Time evolution of $g(t)$.
 % Simulation of $\langle n_a(t_a,t)\rangle$ for the PDF of Eq. (\ref{ageq11}) for slightly aging. The index of $\beta$ for different $\lambda$,   drawn for $\lambda=0. 1,
%0. 05, 0. 03, 0. 02$ and $10^{-10}$.  the parameter $\beta$ ranges from $0. 6$ to $1$. Normal diffusion for $\beta=1$(larger$\lambda$),   and subdiffusion for $\beta=0. 6$ (smaller $\lambda$).  The parameters are $\alpha=0. 6$,  $t_a=3$,  $t=500$.
}\label{agfig91}
\end{figure}

In the analysis of this paper, we use the following function several times.
\begin{equation}\label{age917}
g(z)=z^{\alpha-1}\exp(-\lambda z)E_{\alpha, \alpha}(\lambda^{\alpha}z^{\alpha}).
\end{equation}
When $z \rightarrow +\infty$, from (\ref{ageq915}), there exists
\begin{equation}\label{age918}
g(z) \sim \frac{\lambda^{1-\alpha}}{\alpha};
\end{equation}
see Fig. \ref{agfig91}.

While for $z \rightarrow 0$, from (\ref{age913}), we have
$g(z) \sim 1/\Gamma(\alpha)z^{\alpha-1}$.

\section{Laplace transform of $\varphi(t)$} \label{AppendixB}
Here we present the Laplace transform of Eq. (\ref{ageq11}). From the definition of the survival probability on a site, i.e., the probability that the waiting time on a site exceeds $t$,
\begin{equation}\label{ageq921}
\Psi(t) = \int_{t}^{\infty} \varphi(\tau)d\tau=1-\int_{0}^{t} \varphi(\tau)d\tau.
\end{equation}

Using the definition of incomplete Gamma function, for Eq. (\ref{ageq11}) we have
\begin{equation}\label{ageq922}
\begin{split}
\Psi(t) &\sim \lambda^\alpha\int_{\lambda t}^{\infty} \exp(-z)z^{-\alpha-1} dz \\
       & =\lambda^\alpha\Gamma(-\alpha, \lambda t).
\end{split}
\end{equation}
According to the Laplace transform of the incomplete Gamma function $\mathcal{L}[\Gamma(-\alpha,\lambda t)]= \Gamma(-\alpha)[1-(\frac{u+\lambda}{\lambda})^\alpha]/u$ with $\Re e(-\alpha)>-1$ and  $\Psi(t)$, we have
\begin{equation}\label{ageq923}
\psi(u)  \sim 1+\lambda^\alpha-(u+\lambda)^{\alpha}.
\end{equation}
From Eq. (\ref{ageq923}), we have two useful asymptotics.
For the long time scale, $t \gg 1/\lambda$, (i.e., $u \ll \lambda$) by the Taylor expansion, we have %using the Taylor expand of $(u+\lambda)^\alpha$ in $u$ to second order we have,
\begin{widetext}
\begin{equation}\label{ageq924}
\psi(u)\sim 1-\alpha\lambda^{\alpha-1}u+\alpha(1-\alpha)\lambda^{\alpha-2}u^2+\cdots +(-1)^{n+1}(-\alpha)(1-\alpha)\cdots ((n-1)-\alpha)\lambda^{\alpha-n}u^n.
\end{equation}

Notice that $\psi(0)=1$, so the PDF is normalized. From the definition of $\langle\tau^{n} \rangle=\int_{0}^{\infty} \tau^{n} \psi(\tau)d\tau$, we can get $\langle \tau \rangle=\alpha\lambda^{\alpha-1}$, and $\langle \tau^2 \rangle=(1-\alpha)(-\alpha)\lambda^{\alpha-2}$. For the general cases, $\langle \tau^n \rangle=-\Gamma(n-\alpha)/\Gamma(-\alpha)\lambda^{\alpha-n}$, i.e., $\psi(u) \sim 1-\langle \tau \rangle u+\langle \tau^2 \rangle u^2+\cdots+(-1)^n\langle \tau^n \rangle u^n$.
For short time scale, $t_{0}<<t \ll 1/\lambda$  (i.e. $u \gg \lambda$),, we have
\begin{equation}\label{ageq925}
\psi(u)  \sim 1+u^\alpha\left[\left(\frac{\lambda}{u}\right)^\alpha-\left(1+\frac{\lambda}{u}\right)^\alpha \right]
       \sim 1-u^{\alpha}.
%  \begin{split}
%    \psi(u) & =1+u^\alpha\left[\left(\frac{\lambda}{u}\right)^\alpha-\left(1+\frac{\lambda}{u}\right)^\alpha \right] \\
%      & \sim 1-u^{\alpha}.
%  \end{split}
\end{equation}
\end{widetext}
\end{appendices}
\newpage %Just because of unusual number of tables stacked at end


\begin{thebibliography}{50}
\bibitem{levy:1} The one sided stable L\'{e}vy distribution in Laplace space can be conveniently expressed as, $\int_{0}^{\infty}\exp(-ut)L_{\alpha}(t)dt=\exp(-u^{\alpha})$, for $u>0$.
\bibitem{note:1} The $p_{th} $ moment of a probability distribution of  $n_{a}(t_{a},t)$ with density $P_{n_{a}}(t_{a}, t)$ is defined as, $n_{a}^{p}(t_{a}, t)=\int_{0}^{\infty} n_{a}^{p}(t_{a}, t)P_{n_{a}}(t_{a}, t) d n_{a}. $
\end{thebibliography}
\end{document}